\documentclass[a4paper,11pt]{article}
\usepackage{array}
\usepackage{theorem}
\usepackage{amsmath,amscd,amssymb}
\usepackage{latexsym}
\usepackage{epsfig}
\usepackage{xy}

\theorembodyfont{\sl}

\newtheorem{lemma}{Lemma}[section]

\newtheorem{proposition}[lemma]{Proposition}
\newtheorem{theorem}[lemma]{Theorem}
\newtheorem{corollary}[lemma]{Corollary}

\newcommand{\CC}{\mathbb C}

\newcommand{\HH}{\mathbb H}

\newcommand{\NN}{\mathbb N}

\newcommand{\QQ}{\mathbb Q}
\newcommand{\RR}{\mathbb R}

\newcommand{\ZZ}{\mathbb Z}

\newcommand{\Prod}{\prod\limits}

\newcommand{\Sp}{\mathop{\mathrm {Sp}}\nolimits}

\newcommand{\SL}{\mathop{\mathrm {SL}}\nolimits}

\newcommand{\Ker}{\mathop{\mathrm {Ker}}\nolimits}

\newcommand{\Lift}{\mathop{\mathrm {Lift}}\nolimits}
\newcommand{\Blift}{\mathop{\mathrm {B}}\nolimits}

\renewcommand{\Im}{\mathop{\mathrm {Im}}\nolimits}

\renewcommand{\div}{\mathop{\mathrm {div}}\nolimits}
\newcommand{\Kthree}{\mathop{\mathrm {K3}}\nolimits}
\newcommand{\qedsymbol}{\mbox{$\Box$}}
\newcommand{\qed}{\unskip\nobreak\hfil\penalty50\hskip1em\hbox{}\nobreak
\hfill\qedsymbol\parfillskip=0pt\finalhyphendemerits=0}
\newenvironment{proof}{\begin{ProofwCaption}{Proof}}{\end{ProofwCaption}}
\newenvironment{ProofwCaption}[1]
 {\addvspace\theorempreskipamount \noindent{\it #1.}\rm}
 {\qed \par \addvspace\theorempostskipamount}

\begin{document}

\title{The Siegel modular forms of genus $2$ with the simplest divisor}
\author{V.~Gritsenko and F.~Cl\'ery}
\date{}
\maketitle
\begin{abstract}
We prove that there exist  exactly  eight Siegel modular forms
with respect to the congruence subgroups of Hecke type of the paramodular groups
of genus two vanishing  precisely along the  diagonal of the Siegel upper half-plane.
This is a solution of a question formulated during the conference
``{\it Black holes, Black Rings and Modular Forms}" (ENS, Paris, August 2007).
These modular forms generalize the classical Igusa form and the forms constructed
by Gritsenko and Nikulin in 1998.
\end{abstract}

\bigskip

{\bf{Introduction: dd-modular forms}}
\bigskip

The first cusp form for the Siegel modular group $\Sp_2(\ZZ)$
is the Igusa form $\Psi_{10}$. In fact
$\Psi_{10}=\Delta_5^2$ where $\Delta_5$ is the product
of the ten even theta-constants (see \cite{F}).
This modular form has a lot of remarkable properties.
One of the main features of $\Delta_5$  is that
it vanishes (with order one) precisely along
the  diagonal
$$
{\cal H}_1=
\bigl\{\left(
\begin{array}{cc}
\tau & 0\\ 0 & \omega\end{array}\right),\  \tau,\, \omega \in
\HH_1\bigr\}\subset \HH_2
$$
of the Siegel upper half-plane
$
\HH_2=\{Z=
\left(\smallmatrix
\tau&z\\z&\omega\endsmallmatrix\right)\in M_2(\CC),\, \Im(Z)>0\}.
$
 It is known that
$\Delta_5$ determines a Lorentzian Kac--Moody super Lie algebra
of Borcherds type. See \cite{GN1}--\cite{GN2} where
two lifting constructions of $\Delta_5$ were proposed
$$
\Delta_5(Z)=\Lift(\eta^9(\tau)\vartheta(\tau,z))=\Blift (\phi_{0,1})(Z)
$$
where $\eta$  is  the Dedekind eta-function and  $\vartheta$ is
the Jacobi theta-series (see (\ref{jtheta})).
This relation  gives  the two parts of the denominator identity
for the Borcherds algebra determined by $\Delta_5$.
There exists a geometric interpretation of this identity
in terms of the arithmetic mirror symmetry (see \cite{GN4}).
Moreover $2\phi_{0,1}$ is the elliptic genus of a $\Kthree$ surface and
$\Psi_{10}^{-2}$ is related to the so-called second quantized
elliptic genus of $\Kthree$ surfaces (see \cite{DMVV}, \cite{G4}).
These facts explain the  importance of $\Delta_5$ in the theory
of strings (see \cite{DVV}, \cite{Ka}).
During the conference ``{\it Black holes, Black Rings and Modular Forms}"
(ENS, Paris, August 2007)
there was  formulated a problem
on the existence of Siegel modular forms similar to $\Delta_5$
with respect to the congruence
subgroup of Hecke type
$$
\Gamma^{(2)}_0(N)=\{
\begin{pmatrix}
A&B\\C&D\end{pmatrix}\in \Sp_2(\ZZ)\,|\, \ C\equiv 0\mod N\}.
$$
Such Siegel modular forms can characterize the black holes entropy
and the degeneracy of dyons
for some class of CHL string compactification (see \cite{DJS}, \cite{DS}, \cite{DN}, \cite{DG}).
Mathematically we can reformulate this question  as follows:
{\it  finding all  Siegel modular forms $F$
with respect to  $\Gamma^{(2)}_0(N)$
(with a character or a multiplier system)
such that $F$ vanishes exactly along the $\Gamma^{(2)}_0(N)$-translates
of the diagonal ${\cal H}_1\subset \HH_2$ and with vanishing order one.}
\smallskip

We call such  functions {\bf dd-modular forms}:
 modular forms with the diagonal
divisor.    A dd-modular form is a natural  generalization of  $\Delta_5$.
In this paper we give the complete answer to this  problem.
\begin{theorem}
For the congruence subgroups
$\Gamma_0^{(2)}(N)$ with $N>1$ there are exactly three  dd-modular forms:
$\nabla_3$ of weight $3$ for $\Gamma^{(2)}_0(2)$ with a character of order $2$,
$\nabla_2$ of weight $2$ for $\Gamma^{(2)}_0(3)$ with a character of order $2$
and $\nabla_{3/2}$ of weight $3/2$ for $\Gamma^{(2)}_0(4)$
with a multiplier system of order $4$.
\end{theorem}
In fact we get a result which is stronger than the theorem above.
We give the  full classification of the dd-modular forms
for the  Hecke subgroups $\Gamma_t(N)$ (see (\ref{GtN}))
of the symplectic paramodular groups
$\Gamma_t$. Theorem \ref{dd-class} claims that
there are  exactly  {\it eight} such dd-modular forms.
Four of them,  $\Delta_5$ and the modular forms $\Delta_2$, $\Delta_1$
and $\Delta_{1/2}$ of weights $2$, $1$ and $1/2$ with respect to the paramodular
groups $\Gamma_2$, $\Gamma_3$ and $\Gamma_4$ were constructed  in \cite{GN2}.
The other four modular forms are
the three functions of Theorem 0.1 and the dd-modular form
$Q_1$ of weight $1$  and character of order $4$
with respect to the congruence subgroup
$\Gamma_2(2)$ of the paramodular group $\Gamma_2$.

These eight remarkable modular forms can be considered as the best possible
three dimensional analogues of the Dedekind $\eta$-function.
We expect a number of interesting applications
of these new functions in the string theory, in the theory
of Lorentzian Kac--Moody algebras and in algebraic geometry.

The paper contains three sections.
In \S 1 we prove that there might exist only nine dd-modular forms with respect
to $\Gamma_t(N)$. In \S2 using the lifting of Theorem \ref{J-lift} we construct seven
dd-forms and the square of $\nabla_{3/2}$. Moreover using the particular form
of $Q_1$ we prove that the ninth dd-form does not exist.
In \S 3 using Theorem \ref{product} about the Borcherds automorphic products for congruence subgroup
$\Gamma_{t}(N)$ we construct the last dd-modular form
$\nabla_{3/2}$ of weight $3/2$.

{\it Acknowledgements}:
We would like to thank A. Dabholkar and B. Poline for useful discussions about Siegel modular forms.
The first author is grateful to the  Max-Planck-Institut f\"ur Mathematik
in Bonn for hospitality in 2008 where this work was substantially  done.

\section{Classification of the dd-modular forms}
\label{sec1}

One of the main idea of  our approach (see \cite{G1}--\cite{G3}) is that
in order to understand
better the properties of Siegel modular forms of genus two
one has to consider not only the  modular group $\Gamma_1=\Sp_2(\ZZ)$
and  its congruence subgroups, but the  integral symplectic groups
$\Gamma_t$, the paramodular groups, for all $t\ge 1$.
In this section  we give the complete classification of the dd-modular forms
for the most natural congruence subgroups of the paramodular  groups.
Let $t$ and $N$ be positive integers. We put
\begin{equation}\label{GtN}
\Gamma_t(N)=\biggl\{\begin{pmatrix}
*& *t&  *& *
\\ * & *&*  &*t^{-1}\\
*N&*Nt&*&*\\
*Nt&*Nt&*t&*  \end{pmatrix}\in \Sp_2(\QQ),\quad{\rm all }\  *\in\ZZ \biggr\}.
\end{equation}
The group $\Gamma_t=\Gamma_t(1)$ is conjugated to the integral
symplectic group of the integral skew-symmetric form with elementary
divisors $(1,t)$ (see \cite{G2}, \cite{GH2}).
The quotient
$\Gamma_t\setminus \HH_2$ is the moduli space of the $(1,t)$-polarized
Abelian surfaces.
If $t=1$ then $\Gamma_1=\Sp_2(\ZZ)$ and $\Gamma_1(N)=\Gamma_0^{(2)}(N)$
is the Hecke subgroup from the introduction.

Let $\Gamma_t(N)^+=\Gamma_t(N)\cup \Gamma_t(N)V_t$
be a normal double extension of $\Gamma_t(N)$ in $\Sp_2(\RR)$
where
\begin{equation}\label{Vt}
V_t=\frac{1}{\sqrt{t}}\begin{pmatrix}
0&t&0&0
\\ 1&0&0&0\\
0&0&0&1\\
0&0&t&0  \end{pmatrix}\in \Sp_2(\RR).
\end{equation}
Repeating  the proof of Lemma 2.2 of \cite{G2}  we obtain

\begin{lemma}\label{Gamma+} The group $\Gamma_t(N)^+$ is generated
by  $V_t$ and by its parabolic  subgroup
\begin{equation}\label{Gamma-infty}
\Gamma_{t}^{\infty}(N)=\biggl\{\pm \begin{pmatrix}
* &0&*&*\\
* &1&*&*/t\\
N*&0&*&*\\
0 &0&0&1  \end{pmatrix}\in \Gamma_t(N), \ {\rm all\  }\ \  *\in \ZZ \biggr\}.
\end{equation}
\end{lemma}

Let $\Gamma<\Sp_2(\RR)$ be an arithmetic subgroup.
In this paper $\Gamma$ will be one of the groups $\Gamma_t(N)$ or $\Gamma_t(N)^+$.
{\it A modular form} of weight $k$
($k$ is integral or half-integral) for  the subgroup
$\Gamma$  with a character (or a multiplier system)
$\chi :  \Gamma\to \CC^\times$ is a holomorphic function on $\HH_2$
which satisfies the functional equation
$$
(F|_k \gamma)(Z)=\chi(\gamma)F(Z)\qquad \text{for any }\
\gamma\in \Gamma.
$$
We denote by
${}|_k$ the standard slash operator on the space
of functions on $\HH_2$:
$$
(F|_k\gamma)(Z):=\det(CZ+D)^{-k}F(M\langle Z\rangle)
$$
where
$\gamma=\left(\smallmatrix A&B\\C&D \endsmallmatrix\right)\in \Sp_2(\RR)$
and $M\langle Z\rangle  =(AZ+B)(CZ+D)^{-1}$.
For a half-integral $k$ we choose one of the holomorphic square roots by the
condition
$\sqrt{\det(Z/i)}>0$ for $Z=iY\in \HH_2$.
We denote by $M_k(\Gamma,\chi)$ the finite dimensional  space
of all modular forms of this type.

A dd-modular form is  a particular case of  reflective modular
forms whose divisor is defined by reflections in the corresponding integral
orthogonal group. In \cite{GN2}--\cite{GN3} we classified the reflective modular
forms for the paramodular groups $\Gamma_t$.
In particular we found all Siegel modular forms
for the paramodular groups $\Gamma_t$
with the diagonal divisor. To classify all possible
dd-modular forms for the congruence subgroups  $\Gamma_t(N)$
we use  the  method of multiplicative symmetrization (see \cite[\S 3.1]{GN2}).
The next proposition is a generalization of  Proposition 1.1 of \cite{GH1}
in which we studied the case $N=1$.

\begin{proposition}\label{prop-class}
If $F_k$ is a dd-modular form  of integral
(or half-integral) weight $k$ with a character (or a multiplier system)
with respect to $\Gamma_t(N)$
then the triplet $(t,N;k)$ can take one of the nine values
$$
(1,1;5),\ (2, 1;2),\ (3,1;1),\  (4,1;\frac{1}{2}),\  (1,2;3),
\  (1,3;2),\  (1,4;\frac{3}{2}),
$$
$$
(2,2;1),\  (2,4;\frac{1}{2}).
$$
The corresponding
dd-modular forms are, if they exist, unique up to a scalar.
\end{proposition}
\begin{proof}
Uniqueness of a dd-modular form for a fixed group
follows from the  Koecher principle (see \cite{F}).
Let $F$ be a non-zero modular form of weight $k$ with respect to $\Gamma_t(N)$.
We use the following operator of multiplicative symmetrization
$$
[F]_{1}=
\prod\limits_{\gamma\in
\Gamma^{(int)}_t(N)\setminus \Sp_2(\ZZ)}
F|_k \gamma\quad{\ \ \rm where\ \ } \Gamma^{(int)}_t(N)=\Gamma_t(N)\cap \Sp_2(\ZZ).
$$
It is clear that $[F]_{1}$ is a non-zero modular form
with respect to  $\Gamma_1=\Sp_2(\ZZ)$.
\begin{lemma}\label{ind}
For any integral $t\ge 1$ and $N\ge 1$ we have
$$
[\Gamma_1:\Gamma^{(int)}_{t}(N)]=
\biggl((Nt)^3\prod_{p\,|(tN)}(1+p^{-1})(1+p^{-2})
\biggr)
\cdot \prod_{p\,|(t,N)}(1+p^{-1}).
$$
\end{lemma}
\begin{proof}
The  diagram of the subgroups shows that
$$
[\Gamma_1:\Gamma^{(int)}_{t}(N)]=
\frac{[\Gamma_1:\Gamma^{(int)}_{tN}]\cdot
[\Gamma^{(int)}_{tN} : \Gamma^{(int)}_{tN}\cap\Gamma_t(N)]}
{[\Gamma^{(int)}_{t}(N):\Gamma^{(int)}_{tN}\cap\Gamma_t(N)]}
$$
where  $\Gamma^{(int)}_d=\Gamma_d\cap \Gamma_1$ is a subgroup
of the paramodular group.
It is known (see \cite[\S 1]{GH1}) that
$$
[\Gamma_1:\Gamma^{(int)}_{tN}]=
(tN)^3\prod_{p\,|(tN)}(1+\frac{1}p)(1+\frac{1}{p^2}).
$$
Analyzing the form of the elements in the subgroups we obtain
$$
[\Gamma^{(int)}_{tN}:\Gamma^{(int)}_{tN}\cap\Gamma_t(N)]=
[\SL_2(\ZZ):\Gamma_0(N)]=
N\prod_{p\,|N}\bigl(1+\frac{1}p\bigr),
$$
and
$$
[\Gamma^{(int)}_{t}(N):\Gamma^{(int)}_{t}\cap\Gamma_{tN}]=
[\Gamma_0(t):\Gamma_0(tN)].
$$
This  gives us the formula of the lemma.
\end{proof}

Let
$\pi_{t,N}:\HH_2\rightarrow {\cal A}_t(N)=\Gamma_t(N)\backslash\HH_2$
be the quotient map.
Note that
${\cal A}_{t}={\cal A}_t(1)$ is the moduli space of $(1,t)$-polarized
Abelian surfaces.
For $N=1$ the  image $\pi_{t,1}({\cal H}_1)$ in ${\cal A}_t$
parameterizes split polarized Abelian surfaces.
For $t=N=1$ this is  the Humbert surface $H_1$
of discriminant $1$ in ${\cal A}_{1}$  and one
can consider the divisor $\pi_{1,1}({\cal H}_1)$
as  the discriminant of the moduli space of curves of genus $2$.

Let us assume that $F$ has a diagonal divisor of multiplicity $m\ge 1$,
i.e.,
$\div_{{\cal A}_{t}(N)} F=m\cdot \pi_{t,N}({\cal H}_1)$.
We note that $H_1$ is  irreducible in ${\cal A}_{1}$
(for the theory of Humbert surfaces see \cite{vdG} and \cite{GH2}).
It follows that
$
\div_{{\cal A}_{1}} ([F]_1)
$
is the Humbert surface $H_1$  with some multiplicity $d$.
Therefore according to the Koecher principle
$$
[F]_1(Z)=C\cdot\Delta_5(Z)^d
$$
where $C$ is a non-zero constant.
In order to calculate the multiplicity $d$ we note that
the stabilizer of ${\cal H}_1$ in $\Sp_2(\RR)$
is the group generated by the direct product of two copies of
$\SL_2(\RR)$ in $\Sp_2(\RR)$ and the involution $J$:
\begin{equation}\label{SL2}
\SL_2(\RR)\times \SL_2(\RR) \cong
\biggr\{\left(
\smallmatrix
a&0  &b&0\\
0&a_1&0&b_1\\
c&0  &d&0\\
0&c_1&0&d_1
\endsmallmatrix\right)\in \Sp(\RR)\biggl\},\quad
J=
\left(\smallmatrix
0&1&0&0\\
1&0&0&0\\
0&0&0&1\\
0&0&1&0
\endsmallmatrix\right).
\end{equation}
The order   of zero of $[F]_1$ along $H_1$
is equal to the number of left cosets $\Gamma^{(int)}_t(N)M$
in $\Sp_2(\ZZ)$ containing an element from the stabilizer of ${\cal H}_1$.
Therefore we have to find the number of the distinct
cosets $\Gamma^{(int)}_t(N)M$ with $M\in {\rm St}_{\Sp_2(\ZZ)}({\cal H}_1)$.
The involution $J$ permuting the diagonal elements $\tau$ and $\omega$
of $Z\in \HH_2$ belongs to $\Gamma^{(int)}_t(N)$
if and only if  $t=1$.
If $t>1$ then $\Gamma^{(int)}_t(N)M_1\ne \Gamma^{(int)}_t(N)JM_2$
for any $M_1$, $M_2$ in $\SL_2(\ZZ)\times \SL_2(\ZZ)$.
It gives us the factor $2$ if $t>1$.
Therefore $[F]_1$ vanishes along $H_1$ with order
$$
d=2^{\delta(t)}\,m[\SL_2(\ZZ):\Gamma_0(N)]\cdot[\SL_2(\ZZ):\Gamma_0(tN)].
$$
where $\delta(t)=0$ if $t=1$ and $\delta(t)=1$ if $t>1$.

The weight of the symmetrization$[F]_1$ equals the weight of $F$ multiplied
by the index of the subgroup calculated in the  above lemma.
The relation between $[F]_1$ and $\Delta_5$  gives us the following identity
between  the weights of these modular forms
$$
\biggl(kN\prod_{\substack{p\,|N\vspace{0.5\jot}\\ p\,\nmid t}}
\frac{p^2+1}{p(p+1)}\biggr)\cdot
t^2\prod_{p\,|t}\frac{p^2+1}{p^2}= 2^{\delta(t)}\,5m
$$
where $k\in \ZZ/2$ is the weight of $F$.
For any fixed $m$ simple arguments of divisibility show
that there exist only a few possibilities for $(t,N;k)$.
If $F$ is a dd-modular form of weight $k$ (i.e., if $m=1$),
then there are only four triplets  $(t,N;k)$ with  $t=1$
and five more $(t,N;k)$ for $t>1$.
This proves the proposition.
\end{proof}

In what follows we construct eight dd-modular forms and we prove
that a dd-modular form of type $(2,4;\frac 1{2})$
does not exist.
\begin{theorem}\label{dd-class}
For the  Hecke  congruence subgroups
$\Gamma_t(N)<\Gamma_t$  there are exactly eight dd-modular forms.
They belong to the spaces
$$
M_5(\Gamma_1,\chi_2);\quad
M_2(\Gamma_2,\chi_4),\  M_3(\Gamma_0^{(2)}(2),\chi_2);\quad
M_1(\Gamma_3,\chi_6),\  M_2(\Gamma_0^{(2)}(3),\chi_2);
$$
$$
M_{\frac{1}2}(\Gamma_4,\chi_8), \ M_{\frac{3}2}(\Gamma_0^{(2)}(4),\chi_4);
\qquad M_1(\Gamma_2(2),\chi_4)
$$
where $\chi_d$ is a character (or a multiplier system) of order $d$
of the corresponding modular group.
\end{theorem}
To prove this theorem we describe two lifting constructions for the congruence
subgroups of Hecke type  of the paramodular groups of genus $2$.

\section{Additive construction of dd-modular forms}

The four dd-modular forms for $N=1$ are the modular forms
$\Delta_5$ for $\Sp_2(\ZZ)$ with a character of order $2$,
$\Delta_2$ for $\Gamma_2$ with a character of order $4$,
$\Delta_1$ for $\Gamma_3$ with a character of order $6$
 and $\Delta_{1/2}$ for $\Gamma_4$ with a multiplier system  of order $8$
(see \cite{GN2}).
In this section we construct the dd-modular forms for $N>1$.
For this aim  we use special Jacobi modular forms of index $1/2$
with respect to the Jacobi group of level $N$
$$
\Gamma^J(N)=(\Gamma^{\infty}_t(N)\cap \Sp_2(\ZZ))/\{\pm 1_4\}\cong
\Gamma_0(N)\ltimes H(\ZZ)
$$
(see (\ref{Gamma-infty})).
The Jacobi group is the semi-direct product of
the Heisenberg group
$$
H(\ZZ)=
\biggl\{[\lambda, \mu; \kappa]=
\begin{pmatrix}
1& 0&  0& \mu\
\\ \lambda & 1&\mu  &\kappa\\
0&0&1&-\lambda\\
0&0&0&1
\end{pmatrix}, \quad \lambda, \mu, \kappa\in \ZZ\biggr\}
$$
and
$$
\Gamma_0(N)=\bigl\{
\begin{pmatrix} a&b\\c&d\end{pmatrix}\in \SL_2(\ZZ),\quad c\equiv 0\mod N
\bigr\}.
$$
We embed $\gamma\in \Gamma_0(N)$  in $\Gamma_t(N)$
using the first copy of $\SL_2$ in (\ref{SL2})
(the second term is the unit matrix). We denote this matrix by $\widetilde \gamma$
and we identify $\Gamma_0(N)$ with this subgroup of $\Gamma_t(N)$.

Let $t$ and $k$ be integral or half-integral positive numbers.
A holomorphic function $\phi$ on $\mathbb{H}_1\times\CC$
is called  { \it holomorphic Jacobi form for $\Gamma_0(N)$
of weight $k$ and index $t$} with a character
(or a multiplier system)  $v : \Gamma^J(N) \rightarrow \CC^*$  if the function
$\widetilde\phi(Z):=\phi(\tau,z)e^{2i\pi t\omega}$ of  $Z\in\mathbb{H}_2$
is a $\Gamma^J(N)$-modular form with  character (or multiplier system) $v$,
i.e., if it satisfies
\begin{equation}\label{J-slash}
(\widetilde\phi\big\vert_k \gamma)(Z
)=v(\gamma)\widetilde\phi(Z) \ \  {\text{for any}}\ \
\gamma\in\Gamma^J(N)
\end{equation}
and for each $M\in \SL(2,\ZZ)$ it has a Fourier expansion of the following type
\begin{equation}\label{Fexp}
\bigl(\widetilde{\phi}\big\vert_{k}\widetilde{M}\bigr)(Z)=
\sum_{\substack {n,\, l\vspace{0.5\jot} \\ 4nt-l^2\geqslant 0}}
c_M(n,l)\,q^nr^l s^t
\end{equation}
where $n,l$ are in $\QQ$,
$q=e^{2i\pi\tau}$, $r=e^{2i\pi z}$ and  $s=e^{2i\pi\omega}$.
The last condition means that $\phi$ is holomorphic
at the cusp determined  by $M$ (see \cite{EZ}).
The form $\phi$ is called {\it cusp form} if $c_M(n,l)\ne 0$ only for
$4nt-l^2>0$ for all $M$.
We call the form $\phi$ a {\it weak} Jacobi form if in its Fourier expansions
$c_M(n,l)\neq 0$ only for $n\geqslant 0$.
The Jacobi form  $\phi$ is called {\it nearly holomorphic}
if there exists $n\in\NN$ such that $\Delta^n\phi$ is
a weak Jacobi form
where $\Delta$ is the Ramanujan $\Delta$--function.

We denote by $J_{k,t}(\Gamma_0(N), v)$
the space of all Jacobi forms with a character
(or a   multiplier system)
for
$\Gamma^J(N)=\Gamma_0(N)\ltimes H(\ZZ)$.
We denote the space of corresponding weak (resp. nearly holomorphic)
Jacobi forms by
$J_{k,t}^{w}(\Gamma_0(N))$ (resp. $J_{k,t}^{nh}(\Gamma_0(N))$).
In \S 4, we use nearly holomorphic  Jacobi forms
of weight $0$ for $\Gamma_0(N)$
in order to construct Borcherds products. In this section we work
with holomorphic Jacobi forms.

The main  example of Jacobi forms of half-integral index is
the Jacobi theta-function of level $2$ (see (\ref{thetacar})):
\begin{multline}\label{jtheta}
\vartheta(\tau ,z)=-i\vartheta_{1,1}^{(2)}(\tau,z)=
\sum_{m\in \ZZ}\,\biggl(\frac{-4}{m}\biggr)\, q^{{m^2}/8}\,r^{{m}/2}\\
=-q^{1/8}r^{-1/2}\prod_{n\ge 1}\,(1-q^{n-1} r)(1-q^n r^{-1})(1-q^n)
\end{multline}
is an element of
$J_{\frac{1}{2},\frac{1}{2}}(SL(2,\ZZ),v_{\eta}^3\times v_H)$
where $v_\eta^3$ is the multiplier system of the cube
of the Dedekind $\eta$-function
and
$$
v_H([\lambda, \mu; \kappa])=(-1)^{\lambda+\mu+\lambda\mu+\kappa}.
$$
is a character  of the   Heisenberg group.
(See \cite{GN2}  for more details on the Jacobi forms of half-integral index.)

We denote by  $\chi\times v_H^\varepsilon$
the character of $\Gamma^J(N)$ induced by the character (or multiplier system
of finite order)
$\chi: \Gamma_0(N)\to \CC^*$ and by a power $v^\varepsilon_H:H(\ZZ)\to \{\pm 1\}$.
It is easy to see from the definition that
the non trivial  binary character $v_H$ can appear only if the index $t$
is half-integral.

Let
\begin{equation}\label{phi}
\phi\in J_{k,t}(\Gamma_0(N), \chi\times v_H^{2t})
\end{equation}
where  $k\in \NN$, $t\in \NN/2$, $\chi:\Gamma_0(N)\to \CC^*$
is a character of finite order.
We suppose that
\begin{equation}\label{ker-chi}
\Ker(\chi)\supset \Gamma_1(Nq,q)
\end{equation}
for some $q$ where the last group  is defined as follows $\Gamma_1(Nq,q)=$
$$
\bigl\{\begin{pmatrix}
a&b\\c&d
\end{pmatrix}\in \SL_2(\ZZ),\
c\equiv 0\,{\rm mod}\, Nq,\ b\equiv 0\,{\rm mod}\,q,\
a\equiv d\equiv 1 \,{\rm mod}\,Nq\bigr\}.
$$
For this group we introduce the  Hecke operator
(see \cite[Ch.3]{Sh})
$$
T^{(N)}(m)=\sum_
{\substack{ad=m  \vspace{0.5\jot}
\\(a, Nq)=1 \vspace{0.5\jot}
\\b\,{\rm mod} \,d}}
\Gamma_1(Nq,q)\,\sigma_ a\cdot
\begin{pmatrix} a&qb\\0&d\end{pmatrix}
$$
where $a>0$ and  $\sigma_a\in \SL_2(\ZZ)$ such that
$\sigma_a\equiv
\left(\begin{smallmatrix} a^{-1}&0\\0&a\end{smallmatrix}\right)
\hbox{ mod}\ Nq$.
This element induces the Hecke operator on Jacobi form
$\tilde\phi(Z)=\phi(\tau,z)\exp(2\pi i \omega)$:
\begin{equation}\label{Tm}
\widetilde{\phi}|_k T_-^{(N)}(m)(Z)=
m^{k-1}\sum_
{\substack{ad=m \vspace{0.5\jot}
\\(a, Nq)=1 \vspace{0.5\jot}
\\b\,{\rm mod} \,d}}
d^{-k}\chi(\sigma_a)\phi(\frac{a\tau +bq}d,\,az)
\,e^{2\pi i mt\omega}
\end{equation}
(compare with \cite{EZ} and \cite[(1.11)--(1.12)]{GN2}).

\begin{lemma}\label{T(m)}
Let $\phi$ be as in (\ref{phi}) and (\ref{ker-chi}).
We suppose that $m$ is coprime to $q$ and  that $m$ is odd
if $t$ is  half-integral. Then
$$
\phi|_k\, T_-^{(N)}(m)\in J_{k,mt}(\Gamma_0(N),\chi_m\times v_H^{2t})
$$
where  $\chi_m$ is a character of $\SL_2(\ZZ)$ defined by
$$
\chi_{m}(\alpha):=\chi(\alpha_m).
$$
For any $\alpha=
\left(\begin{smallmatrix} a&b\\c&d\end{smallmatrix}\right)
\in \Gamma_0(N)$ the matrix $\alpha_m\in \Gamma_0(N)$ is defined
by the condition
$$
\alpha_m\equiv
\begin{pmatrix} a\mod Nq&m^{-1}b\mod q\\mc\mod Nq&d \mod Nq\end{pmatrix}.
$$
\end{lemma}
The proof of the lemma is similar to the proof of \cite[Lemma 1.7]{GN2}.
One has  to use that $\Gamma_1(Nq,q)$ is a normal subgroup
of $\Gamma_0(N)$.
We note that we do not assume that $m$ is coprime to $N$.

In this paper we consider  Jacobi forms with  special
characters such that
$\Ker(\chi)\supset \Gamma_1(Nq,q)$.
If $q=1$, then $\chi$   is induced by
a Dirichlet character $\chi_N$ modulo $N$:
\begin{equation}\label{chiN}
\chi\bigr(\left(\smallmatrix a&b\\c&d \endsmallmatrix\right)\bigl)=\chi_N(d).
\end{equation}
To construct all  dd-modular forms
we have to use characters which appear in the theory
of $\eta$-products.
For example, we shall use the following Jacobi forms
$$
\eta(\tau)\eta(2\tau)^4\vartheta(\tau,z)
\quad\text{or}\quad
\frac{\eta(2\tau)^2\eta(4\tau)^4}{\eta(\tau)^2}\,\vartheta(\tau,z)^2
$$
(see the proof of Theorem \ref{dd-class} below).
The corresponding  characters
can be calculated using the conjugation of the  multiplier system $v_\eta$
of order $24$ of the Dedekind eta-function.
This explains the role of  the number $24$
in the  lifting construction of Theorem \ref{J-lift}.
This theorem  generalizes to congruence subgroups
the lifting constructions of  \cite{G1} and \cite{GN2}.

\begin{theorem}\label{J-lift}
Let $\phi\in J_{k,t}(\Gamma_0(N), \chi\times v_H^{2t})$ be a holomorphic Jacobi form
where  $k\in \NN$, $t\in \NN/2$ and  $\chi:\Gamma_0(N)\to \CC^*$
is a character of finite order such that
$\Ker(\chi)\supset \Gamma_1(Nq,q)$.
We assume that $q$ is a divisor of $24$, $qt\in \NN$ and
$\chi\bigr(\left(\smallmatrix 1&1\\0&1\endsmallmatrix\right)\bigl)
=e^{\frac{2\pi i}q}$.

{\bf 1.} Let $q>1$ or $q=1$ and $c(0,0)=0$ where
$c(0,0)$ is the constant coefficient in the Fourier expansion of $\phi$ at $\infty$.
We  fix $\mu\in (\ZZ/q\ZZ)^\times$.
Then the function
$$
F_\phi(Z)={\rm Lift}_\mu(\phi)(Z)
=\sum_{\substack{ m\equiv \mu \, {\rm mod} \,q  \vspace{0.5\jot}\\
m>0}}
\widetilde{\phi}\,|_k\, T_-^{(N)}(m)(Z)
$$
is a modular form for $\Gamma_{qt}(N)^{+}$
with a character $\chi_{t,\mu}$.
The lifting  is a cusp form if $\phi$ is a cusp form.
If $\mu=1$, then ${\rm Lift}(\phi)={\rm Lift}_1(\phi)\not\equiv 0$
for $\phi\not\equiv 0$.
If ${\rm Lift}_\mu(\phi)\not\equiv 0$, then
the character $\chi_{t,\mu}$ is induced by
the character $\chi_\mu\times v_H^{2t}$ of the Jacobi group,
where $\chi_{\mu}$ is a character of $\SL_2(\ZZ)$
$\mu$-conjugated to $\chi$ (see Lemma \ref{T(m)}), and by the relations
$$
\chi_{t,\mu}(V_{qt})=(-1)^k,\quad
\chi_{t,\mu}([0,0;\frac{\kappa}{qt}])=
\exp{(2\pi i\, \frac{\mu\kappa}q)}\quad (\kappa\in \ZZ).
$$
{\bf 2.} Let  $q=1$ and $c(0,0)\ne 0$. We assume that the character $\chi$
of $\Gamma_0(N)$ is induced by a primitive Dirichlet character
$\chi_N$ modulo $N$ (see (\ref{chiN})).
Then
$$
F_\phi(Z)={\rm Lift}(\phi)(Z)
=c(0,0)E_k(\tau,\chi_N)+\sum_{m\ge 1}
\widetilde{\phi}\,|_k\, T_-^{(N)}(m)(Z)
$$
where
$$
E_k(\tau,\chi_N)=2^{-1}L(1-k,\chi_N)+\sum_{n\ge 1} \sum_{a|n} \chi_N(a)\,a^{k-1}
\exp{(2\pi i n\tau)}
$$
is the  Eisenstein series of weight $k$ for $\Gamma_0(N)$
with character $\chi_N$.
\end{theorem}

\noindent
{\bf Remark.}  There is a variant of this theorem if $qt$ is half-integral.
One has to add a conjugation with respect to an element
of the symplectic group over $\QQ$ in order to obtain a modular form
for a congruence subgroup of the paramodular group $\Gamma_{4qt}$.
(See \cite[Theorem 1.12]{GN2}).

\begin{proof} First we prove the convergence of the series
defining $\hbox{Lift}_\mu(\phi)$.
We put
$$
Z=X+iY=\begin{pmatrix}\tau&z\\z&\omega\end{pmatrix}=
\begin{pmatrix}u&x\\x&u_1\end{pmatrix}
+i
\begin{pmatrix}v&y\\y&v_1\end{pmatrix}\in \HH_2.
$$
Then $\det Y=vv_1-y^2=v(v_1-\frac{y^2}v)=v\cdot \tilde v$
where   $\tilde v$ is invariant under the action
of the Jacobi group.
If $\phi$ is a holomorphic function with
Fourier expansion of type (\ref{Fexp}), then
$$
|\phi(\tau ,z)e^{2\pi i t\omega}|e^{2\pi t\tilde v}
= |\phi(\tau ,z)|e^{-2\pi t y^2/v}
$$
does not depend on $u_1$ and $\tilde v$ and
it is bounded
in the domain $v>\varepsilon$ (see \cite{Kl}).
We introduce
$$
\widetilde{\psi}(Z)=\sum_{M_i\in \Gamma_0(N)\setminus \SL_2(\ZZ)}
|\widetilde{\phi}|_k {\widetilde {M_i}}(Z)|.
$$
Then  the function $\widetilde\psi$ is $|_k$-invariant with respect
to the full Jacobi group $\Gamma^J=\SL_2(\ZZ)\ltimes H(\ZZ)$
and $\widetilde{\psi}(Z)e^{2\pi t\tilde v}$
(depending  only on $\tau$ and $z$)
is bounded for  $v>\varepsilon$.
If $0<v<\varepsilon$ then there exists
$M=
\left(\smallmatrix a&b\\c&d \endsmallmatrix\right)
\in \SL_2(\ZZ)$
such that $\Im (M\langle \tau\rangle)>\varepsilon$ and
$$
\widetilde \psi(Z)=|c\tau+d|^{-k} \widetilde \psi(\widetilde M\langle Z\rangle).
$$
Therefore $\widetilde{\psi}(Z)e^{2\pi t\tilde v}={\rm O}(v^{-k})$
if $v\to 0$. Using, if necessary, this estimation
for all  terms $\psi(\frac {a\tau+b}d, az)$
we see that the series which defines  $F_\phi$
has a majorant of type
$$
C\sum_{m\ge 1} m^{k+1} e^{-2\pi tmv_1}
$$
in any compact subset of $\HH_2$.
If $c(0,0)\ne 0$ then we add an Eisenstein series
$E_k(\tau,\chi_N)$ with respect
to $\Gamma_0(N)$ in the lifting construction.
This series is well defined for any weight $k\ge 1$ according to Hecke
(see \cite[Proposition 5.1.2]{Hi}).

Now we prove that the lifting is a modular form.
According to the conditions of the theorem
$\phi$ has the Fourier expansion of the following type
$$
\phi(\tau  ,\,z )
=\sum_{\substack{\, l\equiv 2t \,{\rm mod}\,2
\vspace{0.5\jot}\\n\ge 0,\ n\equiv 1\,{\rm mod}\,q
\vspace{0.5\jot}\\\frac{4nt}q\ge  \frac {l^2}4}}
c(n,l)\,\exp{\bigl(2\pi i (\frac n{q}\tau +\frac l{2}z)\bigr)}.
$$
By the definition (\ref{Tm}) of the Hecke operators
we have
\begin{equation}
\bigl(\widetilde{\phi}|_k \,T_{-}^{(N)}(m)
\bigr)(Z)=
\end{equation}
$$
m^{k-1}\hspace{-0.2cm}\sum_
{\substack{ad=m  \vspace{0.5\jot}
\\(a, Nq)=1 \vspace{0.5\jot}
\\b\,{\rm mod} \,d}}\hspace{-0.2cm}
d^{-k} \chi(\sigma_a)\hspace{-0.3cm}
\sum_{\substack{\, l\equiv 2t \,{\rm mod}\,2
\vspace{0.5\jot}\\n\ge 0,\ n\equiv 1\,{\rm mod}\,q}}
\hspace{-0.3cm} c(n,l)\,\exp{\bigl(2\pi i (\frac {n(a\tau+bq)}{dq}
+\frac {al}{2}z+mt\omega)\bigr)}=
$$
$$
\sum_
{\substack{ad=m\vspace{0.5\jot} \\(a, Nq)=1}}
a^{k-1} \chi(\sigma_a)
\sum_{\substack{\, l\equiv 2t \,{\rm mod}\,2,\ n_1\ge 0
\vspace{0.5\jot}\\\,dn_1\equiv 1\,{\rm mod}\,q}}
c(dn_1,l)\,\exp{\bigl(2\pi i (\frac {an_1}{q}\tau +\frac {al}{2}z
+adt\omega)\bigr)}.
$$
The Jacobi form $\phi$ has a nontrivial character $v_H$
of the Heisenberg subgroup of the Jacobi group if and only if
$2t\equiv 1\mod 2$.
If $t$ is half-integral, then $q$ is pair because
$tq\in \NN$. Therefore in this case for any $m$ coprime to $q$
the character of the Heisenberg group   of the  Jacobi form
$\widetilde\phi|_k T_-^{(N)}(m)$  is equal to $v_H$.
The $\Gamma_0(N)$-part of the character of $\widetilde\phi|_k T_-^{(N)}(m)$
depends only on $m$ modulo $q$  according to Lemma 2.1
because $\phi$  satisfies (\ref{ker-chi}).
In the definition of $F_\phi$ we have   $m\equiv \mu \mod q$.
Therefore, if $m\equiv \mu \mod q$, then
$\widetilde\phi|_k T_-^{(N)}(m)$ is a Jacobi form with character
$\chi_\mu\times v_H^{2t}$.

The number $q$ is a divisor of $24$ and $\mu$ is coprime to $q$.
For any  $x\in (\ZZ/24\ZZ)^\times$ we have $x^2\equiv 1 \mod 24$.
($24$ is the maximal number with this property.
The  same is true for any divisor of $24$.)
Therefore in the formula for the  Fourier expansion of  $\widetilde\phi|_k T_-^{(N)}(m)$
we have that  the coefficient at $\tau$ under the exponent satisfies
the relations $an_1=\frac{an}{d}\equiv mn\equiv \mu \mod q$.

Now we assume that $c(0,0)=0$. Taking the summation over all positive
$m\equiv \mu \mod  q$  we get
$$
F_\phi(Z)=
$$
$$
\sum_
{\substack{a>0,\, d>0\vspace{0.5\jot} \\ad\equiv \mu\, {\rm mod}\,q
\vspace{0.5\jot}\\ (a,N)=1}}\hspace{-0.2cm}
a^{k-1} \chi(\sigma_a)\hspace{-0.3cm}
\sum_{\substack{\, l\equiv 2t \,{\rm mod}\,2
\vspace{0.5\jot}\\n_1>0,\ dn_1\equiv 1\,{\rm mod}\,q}}
c(dn_1,l)\,\exp{\bigl(2\pi i (\frac {an_1}{q}\tau +\frac {al}{2}z
+adt\omega)\bigr)}=
$$
$$
\sum_
{\substack{n,\,m>0
\vspace{0.5\jot}\\
n,m\equiv \mu \, {\rm mod}\,q \vspace{0.5\jot}\\
l\equiv 2t \,{\rm mod}\,2\vspace{1\jot}\\
\frac{4nmt}q \ge \frac{l^2}4}}\hspace{2\jot}
\sum_{\substack{a|(n,l,m)\vspace{0.5\jot}\\ (a,N)=1
\vspace{0.5\jot}\\ a>0}}
a^{k-1}\, \chi(\sigma_a)\,
c(\frac{nm}{a^2}, \frac{l}a)
\,\exp{\bigl(2\pi i (\frac n{q}\tau +\frac {l}{2}z +mt\omega)\bigr)}.
$$
The Jacobi forms
$\widetilde{\phi}|_k\,T_-^{(N)}(m)$ are modular forms with
respect the Jacobi group $\Gamma^J(N)$.
The parabolic subgroup
$\Gamma_{qt}^\infty(N)$ (see (\ref{Gamma-infty}))
differs from
the Jacobi group $\Gamma^J(N)$ by its center.
For $m\equiv \mu \mod q$ the action of the center
is given by
$$
\bigl(\widetilde{\phi}|_k\,T_-^{(N)}(m)\bigl) \big|_k\,[0,0;\frac \kappa{qt}]=
\exp{(2\pi i \frac {\kappa\mu}{q})}\,
\bigl(\widetilde{\phi}|_k\,T_-^{(N)}(m)\bigr).
$$
Therefore the lifting $F_\phi(Z)$ is
a $\Gamma_{qt}^\infty(N)$-modular form of weight $k$
with character
$\chi_\mu\times v_H^{2t}\times \exp{(2\pi i \frac {*\mu}{q})}$.
The Fourier expansion of $F_\phi$ is also invariant
under the transformation
$\{\tau \to qt\omega ,\  \omega \to (qt)^{-1}\tau \}$.
It is induced by   $V_{qt}$
(see (\ref{Vt})). Therefore
$$
(F_{\phi}|_k\, V_{qt}) (Z)=(-1)^kF_{\phi}(Z).
$$
The subgroup $\Gamma_{qt}^\infty(N)$ and $V_{qt}$
generate the group $\Gamma_{qt}^+(N)$ (see Lemma \ref{Gamma+})
and the lifting is a $\Gamma_{qt}^+(N)$-modular form if $c(0,0)\ne 0$.

If $\phi$ is a cusp form then $\phi|_k T_-^{(N)}(m)$ is also a Jacobi cusp form.
In order to prove this we note that $T^{(N)}(m)$ is a part of the full
Hecke operator for the congruence subgroup $\Gamma_0(N)$ (see (\ref{Tmc})).
Therefore for any $M\in \SL_2(\ZZ)$ we have
$
(\widetilde\phi|_k T_-^{(N)}(m))|_k\widetilde M=
\sum_i (\widetilde\phi|_k \widetilde M_i)|_k P_i
$
for some  $M_i\in \SL_2(\ZZ)$ and integral upper triangular matrices  $P_i$
with $\det P_i=m$. All indices of the Fourier coefficients
of $\widetilde\phi|_k \widetilde M_i$ have positive hyperbolic norm
$4nt-l^2>0$ like in (\ref{Fexp}).
The action by upper triangular $P_i$ does not change this property.
It follows that  $\widetilde\phi|_k T_-^{(N)}(m)$
is a cusp form.
We have proved that the index $(n,l)$ of arbitrary    non-zero Fourier coefficient
of the lifting is non degenerate (i.e., $4nt-l^2>0$) for all $0$-dimensional cusps
of the $1$-dimensional cusp determined by  $\Gamma^J(N)$.
$\Gamma_t(N)$ and the full Jacobi group $SL_2(\ZZ)\ltimes H(\ZZ)$
generate the paramodular group $\Gamma_t$.
In order to obtain all cusps of $\Gamma_t(N)$ we can use
the parabolic subgroup
$
\Gamma_\infty=\{P=\left(\smallmatrix A&B\\0&D\endsmallmatrix\right)\in \Sp_2(\ZZ)\}
$
because $\langle\Gamma_t, \Gamma_\infty\rangle=\Sp_2(\ZZ)$ (see \cite{G2}).
We have considered above the action of $\SL_2(\ZZ)$ on the lifting.
After the action by any upper triangular matrix $P$  Fourier coefficients
with degenerate index do not appear. Therefore the lifting of a Jacobi cusp form
is a cusp form.

Let us consider the case  $c(0,0)\ne 0$.
Since $\chi\bigr(\left(\smallmatrix 1&1\\0&1\endsmallmatrix\right)\bigl)
=e^{\frac{2\pi i}q}$  we have that
$q=1$, $t\in \NN$,
the character $\chi$ is induced by a Dirichlet character $\chi_N$ modulo $N$
(see (\ref{chiN})) and $(-1)^k=\chi_N(-1)$.
In the sum
$\sum_{m\ge 1} \widetilde{\phi}\,|_k\, T_-^{(N)}(m)(Z)$ we have
an additional term
$$
c(0,0)\sum_{m\ge 1} \sum_{a|m} \chi_N(a)\,a^{k-1}
\exp{(2\pi i mt\omega)}.
$$
To make the lifting  invariant with respect to  $V_t$
($\omega\mapsto \tau/t$, $\tau\mapsto t\omega$)
we have to add a similar term with respect to  $\tau$. For that we use
the Eisenstein series
$$
E_k(\tau,\chi_N)=2^{-1}L(1-k,\chi_N)+\sum_{n\ge 1} \sum_{a|n} \chi_N(a)\,a^{k-1}
e^{2\pi i n\tau}
\in M_k(\Gamma_0(N),\chi_N)
$$
(see \cite[Proposition 5.1.2]{Hi}).
The Eisenstein series  is a Jacobi
form of weight $k$ and index $0$.
The theorem is proved.
\end{proof}

\noindent
{\bf Proof of Theorem \ref{dd-class}.} We consider the nine possibilities
for dd-modular forms given in Proposition \ref{prop-class}.

{\bf 1}.  $N=1$. The dd-modular forms for the full paramodular group $\Gamma_t$
with $t=1$, $2$, $3$, $4$ were constructed in \cite{GN1}--\cite{GN2}:
$$
\Delta_5(Z)=\Lift(\eta(\tau)^9\vartheta(\tau,z))
\in M_5(\Gamma_1, v_\eta^{12}\times v_H),
$$
\begin{equation}\label{Delta2}
\Delta_2(Z)=\Lift(\eta(\tau)^3\vartheta(\tau,z))\in M_2(\Gamma_2, v_\eta^{6}\times v_H),
\end{equation}
$$
\Delta_1(Z)=\Lift(\eta(\tau)\vartheta(\tau,z))\in M_1(\Gamma_3, v_\eta^{4}\times v_H).
$$
They are cusp forms with character of order $2$, $4$ and $6$ respectively.
Moreover
\begin{equation}\label{Delta1/2}
\Delta_{{1}/2}(Z)={\rm Trivial-}\Lift
(\vartheta(\tau,z))\in M_{{1}/2}(\Gamma_4, v_\eta^{3}\times v_H)
\end{equation}
is the most odd Siegel even theta-function $\theta_{1111}(Z)$ of level $2$
which is a modular form of weight $1/{2}$ and a multiplier system of degree $8$
with respect to $\Gamma_4$.
We construct below the four new Siegel dd-modular forms
$\nabla_3$, $\nabla_2$, $\nabla_{{3}/2}$ and $Q_1$
for the congruence subgroups.
The index denotes the weight of the corresponding modular form.

{\bf 2}. Let $N=2$. Two groups of level $N=2$  appear in Proposition \ref{prop-class}.
We consider two Jacobi forms of index $\frac{1}2$ with respect
to the Hecke congruence subgroup $\Gamma_0(2)$:
$$
\eta(\tau)\eta(2\tau)^4\,\vartheta(\tau,z)\in
J_{3,\frac{1}2}^{cusp}(\Gamma_0(2),\chi^{(2)}_2\times v_H),
$$
$$
\frac{\eta(2\tau)^2}{\eta(\tau)}\,\vartheta(\tau,z)
\in J_{1,\frac{1}2}(\Gamma_0(2),\chi^{(2)}_4\times v_H).
$$
Every cusp $p$ of $\Gamma_0(N)$ has a representative of the form
$p={a}/c$ where $c$ is a positive divisor of $N$ and
$a$ is taken $\mod(c,\frac{N}c)$.
For any divisor $n$ of $N$ the order of $\eta(n\tau)$ at $p$ is equal to
$\frac{(c,n)^2}{24n}$. Using this we check that
the $\Gamma_0(2)$-modular form $\frac{\eta(2\tau)^2}{\eta(\tau)}$
has a zero of order $\frac{1}8$ at $p=\infty$ and is equal to $1/2$ at the second cusp.

The Jacobi theta-series $\vartheta(\tau,z)$ has the multiplier system
$v_\eta^3\times v_H$ of order $8$.
$\eta(2\tau)^8\eta(\tau)^8$ is a well known example of the  modular forms
with respect to $\Gamma_0(2)$. The powers
$\eta(2\tau)^4\eta(\tau)^4$ and $\eta(2\tau)^2\eta(\tau)^2$ are cusp forms
for $\Gamma_0(2)$ with characters $\chi_2^{(2)}$ and $\chi_4^{(2)}$
of $\Gamma_0(2)$ of order $2$ and $4$ respectively.
Using the exact formula for $v_\eta^2$ (see, for example, \cite[Lemma 1.2]{GN2})
we obtain that
$$
\chi_2^{(2)}\bigr(\left(
\smallmatrix a&b\\2c&d\endsmallmatrix\right)\bigl)=(-1)^{b-c},
\qquad
\chi_4^{(2)}
\bigr(\left(\smallmatrix a&b\\2c&d\endsmallmatrix\right)\bigl)
=e^{\frac{2\pi i}4 d(b-c)}
$$
for any matrix in $\Gamma_0(2)$.
In particular
$$
\ker \chi_2^{(2)}\supset \Gamma_1(4,2),\quad
\ker \chi_4^{(2)}\supset \Gamma_1(8,4),
$$
$$
\chi_2^{(2)}\bigr(\left(\smallmatrix 1&1\\0&1\endsmallmatrix\right)\bigl)
=e^{\frac{2\pi i}2}, \quad
\chi_4^{(2)}\bigr(\left(\smallmatrix 1&1\\0&1\endsmallmatrix\right)\bigl)
=e^{\frac{2\pi i}4}.
$$
The lifting construction gives us
\begin{equation}\label{nabla3}
\nabla_3:=\Lift(\eta(\tau)\eta(2\tau)^4\,\vartheta(\tau,z))\in
S_3(\Gamma_0^{(2)}(2), \chi_2^{(2)}\times v_H),
\end{equation}
\begin{equation}\label{Q1}
Q_1:=\Lift(\frac{\eta(2\tau)^2}{\eta(\tau)}\,\vartheta(\tau,z))
\in M_1(\Gamma_2(2), \chi_4^{(2)}\times v_H).
\end{equation}

{\bf 3}. Let $N=3$.
It is known that
$\eta(3\tau)^6\eta(\tau)^6\in S_6(\Gamma_0(3))$.
We consider
$$
\eta(3\tau)^3\vartheta(\tau,z)
\in J_{2,\frac{1}2}^{cusp}(\Gamma_0(3),\chi^{(3)}_2\times v_H)
$$
where $\chi^{(3)}_2$ is a character of order $2$.
Similar to the case $N=2$ one can check that
$$
\chi^{(3)}_2(M)
=(-1)^{a+d+1}\left(\frac{d}{3}\right)\quad
M=\begin{pmatrix} a&b\\3c&d\end{pmatrix}\in \Gamma_0(3)\quad\text{if } \  c\equiv 1 \mod 2
$$
and
$$
\chi^{(3)}_2(M)
=(-1)^{b}\left(\frac{d}{3}\right)\quad\text{if }\  c\equiv 0 \mod 2.
$$
Therefore
$$
\ker \chi_2^{(3)}\supset \Gamma_1(6,2),\quad
\chi_2^{(3)}\bigr(\left(\begin{smallmatrix}
 1&1\\0&1\end{smallmatrix}\right)\bigl)
=e^{\pi i}.
$$
Applying Theorem \ref{J-lift} we obtain
\begin{equation}\label{nabla2}
\nabla_2:=\Lift(\eta(3\tau)^3\vartheta(\tau,z))\in S_2(\Gamma_0^{(2)}(3), \chi_2^{(3)}\times v_H).
\end{equation}

{\bf 4}. Let  $N=4$. We define
\begin{equation}\label{h3/2}
h_{\frac{3}2}(\tau,z)=
\frac{\eta(2\tau)\eta(4\tau)^2}{\eta(\tau)}\,
\vartheta(\tau,z) \in J_{\frac{3}2,\frac{1}2}(\Gamma_0(4),\chi^{(4)}_4\times v_H)
\end{equation}
where $\chi^{(4)}_4$ is a multiplier system of order $4$.
The $\Gamma_0(4)$-modular form
$\eta(2\tau)\eta(4\tau)^2/ \eta(\tau)$
vanishes at two cusps $\infty$ and $\frac{1}2$ and takes non-zero value
at $1$.
Let us consider
$$
h_{\frac{3}2}(\tau,z)^2=
\frac{\eta(2\tau)^2\eta(4\tau)^4}{\eta(\tau)^2}\,
\vartheta(\tau,z)^2 \in J_{3,1}(\Gamma_0(4),\chi^{(4)}_2).
$$
The $\Gamma_0(4)$-modular form
$\eta(4\tau)^4\eta(2\tau)^2\eta(\tau)^4$ has the non-trivial quadratic
character modulo $4$
$$
\chi_2^{(4)}(\left(\smallmatrix a&b\\4c&d\endsmallmatrix\right))
=(-1)^{\frac{(d-1)}{2}}
\quad\text{ where }\ \
\left(\smallmatrix a&b\\4c&d\endsmallmatrix\right)\in \Gamma_0(4).
$$
We have
\begin{equation}\label{F3}
F_3:=\Lift(h_{{3}/2}^2)\in M_3(\Gamma_0^{(2)}(4), \chi_2^{(4)}).
\end{equation}
The Jacobi form $h_{{3}/2}(\tau,z)^2$ has zero of order $2$ for $z=0$.
The Hecke operators of the lifting keep this divisor.
Therefore $F_3$ vanishes with order $2$ along  ${\cal H}_1=\{z=0\}$.
The proof of Proposition \ref{prop-class} shows that
$$
\div_{\HH_2} (F_3)=2\bigl(\bigcup _{\gamma\in \Gamma_0^{(2)}(4)}
\gamma\langle {\cal H}_1\rangle\bigr).
$$
In the next section
we construct a modular form $\nabla_{{3}/2}$ such that
$F_3=\nabla_{3/{2}}^2$ using the Borcherds automorphic products.

{\bf 5}. The last case of  Proposition \ref{prop-class} is  a possible  dd-modular form
of type $(N,t;k)=(4,2;\frac{1}2)$.
Using the exact construction of the dd-modular form  $Q_1$
for $\Gamma_2(2)$
we   prove that a dd-modular form  of weight $\frac{1}2$ with respect to  $\Gamma_2(4)$
does not exist.

Let assume that $D\in  M_{\frac 1{2}}(\Gamma_2(4),\chi)$ is a dd-modular form.
$Q_1$ can be considered as  a modular form with respect to $\Gamma_2(4)<\Gamma_2(2)$.
$Q_1$ vanishes along  the diagonal ${\cal H}_1$ but
its divisor modulo $\Gamma_2(4)$ contains several irreducible components.
Therefore $F=Q_1/D$ is a holomorphic function on $\HH_2$ and it is a $\Gamma_2(4)$-modular
form of weight $\frac 1{2}$ according to the Koecher principle.
Then it has the following Fourier-Jacobi expansion
$$
F(Z)=f_0(\tau)+\sum_{m\ge \frac{1}2} f_{\frac 1{2},m}(\tau,z)\exp(2\pi i m \omega)
$$
where the constant term  $f_0$
is a modular form of weight $\frac{1}2$ with respect to $\Gamma_0(4)$.
The zeroth Fourier-Jacobi coefficient
of D is identically equal to zero
$$
d_0(\tau)=\lim_{v_1\to \infty}
D\bigr(\begin{pmatrix} \tau&z\\z&iv_1\end{pmatrix}\bigl)
\equiv 0
$$
because $D$ is zero for $z=0$.
Considering the Fourier-Jacobi expansions of the both part of the identity
$Q_1=D\cdot F$ we obtain that
$$
\frac{\eta(2\tau)^2}{\eta(\tau)}\,\vartheta(\tau,z)
=f_0(\tau)\cdot d_{\frac 1{2},\frac{1}2}(\tau,z).
$$
Therefore the first non-trivial Fourier-Jacobi coefficient
$d_{\frac 1{2},\frac{1}2}$ of $D$
is equal to
$g\vartheta$ where
$g$ is an automorphic form of weight $0$ with respect to $\Gamma_0(4)$.
The Jacobi theta-series is a modular form of singular weight $1/2$.
For every  Fourier coefficient $c(n,l)$ of $\vartheta$
we have $2n^2-l^2=0$ (see (\ref{jtheta})).
The automorphic form $g$ has a pole at  some cusp. Therefore
the Jacobi form $g\vartheta$ cannot be holomorphic at this cusp.
It follows that  $D$ is not holomorphic.
We finish the proof of Theorem \ref{dd-class} modulo existence
of the dd-modular form $\nabla_{3/2}$.
\smallskip

We note that two new dd-modular forms  $\nabla_2$ and $Q_1$
have  elementary formulae for the Fourier coefficients.
You can compare them with cusp forms $\Delta_2\in S_2(\Gamma_2,\chi_4)$
and $\Delta_1\in S_1(\Gamma_3,\chi_6)$ (see \cite{GN1} and
\cite[Example 1.14]{GN2}).
According to Euler and Jacobi
$$
\eta(\tau)^3=\sum_{n>0}\left(\frac{-4}n\right)nq^{n^2/8}.
$$
Then we obtain
\begin{equation}\label{nabla2-F}
\nabla_2(Z)=\sum_{N>0}\sum_{\substack{m,\,n\in 2\NN+1\vspace{0.5\jot}\\
 3N^2=4mn-l^2}}
N\left(\frac{-4}{Nl}\right)
\sum_{\substack{a|(l,m,n) \vspace{0.5\jot}\\ a>0}}
a\left(\frac{a}3\right)q^{\frac{n}{2}}r^{\frac{l}{2}}s^{\frac{m}{2}}
\end{equation}
because in the lifting formula $\chi(\sigma_a)=\left(\frac{a}3\right)$.
To calculate $Q_1$ we note that
$\eta(2\tau)^2/\eta(\tau)=\frac{1}{2}\vartheta_{1,0}^{(2)}(\tau,0)$
where
$$
\vartheta_{1,0}^{(2)}(\tau,0)=
q^{\frac{1}{8}}\prod_{n\ge 1}(1-q^n)(1+q^{n-1})(1+q^n)
=2\sum_{n\in \NN}q^{\frac{(2n+1)^2}{8}}.
$$
Using the last formula we obtain that
\begin{equation}\label{Q1-F}
Q_1(Z)=\sum_{N>0}
\sum_{\substack{n,\,m \in 4\NN+1\vspace{0.5\jot}
\\ l \in 2\ZZ+1\vspace{0.5\jot}\\(2N+1)^2=2mn-l^2}}
\left(\frac{-4}l\right)
\sigma_0((n,l,m))\,
q^{\frac{n}{4}}r^{\frac{l}{2}}s^{\frac{m}{2}}
\end{equation}
where $\sigma_0((n,l,m))$ is the number of divisors of the greatest common divisor
of $n$, $l$, $m$.

From the proof of theorem given above  we obtain also a description
of the squares of dd-forms as liftings.
\begin{corollary}\label{dd-powers}
The following identities are true
$$
\nabla_3(Z)^2=\Lift\bigl(
\eta(\tau)^2\eta(2\tau)^8\,\vartheta(\tau,z)^2
\bigr)\in S_6(\Gamma_0(2)),
$$
$$
\nabla_2(Z)^2=\Lift\bigl(
\eta(3\tau)^6\vartheta(\tau,z)^2\bigl)\in S_4(\Gamma_0(3)),
$$
$$
Q_1(Z)^2=\Lift\bigl(
\frac{\eta(2\tau)^4}{\eta(\tau)^2}\,\vartheta(\tau,z)^2\bigr)
\in M_2(\Gamma_2(2),\chi_2),
$$
$$
Q_1(Z)^4=\Lift\bigl(
\frac{\eta(2\tau)^8}{\eta(\tau)^4}\,\vartheta(\tau,z)^4\bigr)
\in M_4(\Gamma_2(2)).
$$
\end{corollary}
\begin{proof}
All identities are similar. We prove the last one.
First  $f_1^4\in J_{4,2}(\Gamma_0(2))$.
According to Theorem \ref{J-lift}
we get  $\Lift(f_1^4)\in M_4(\Gamma_2(2))$.
The Jacobi form $f_1^4$ has zero of order $4$ along $z=0$.
The Hecke operators in the lifting construction preserve this divisor.
Therefore the quotient $\Lift(f_1^4)/Q_1^4$ is a constant according to
the Koecher principle. This constant is one. To see this we compare the
first Fourier-Jacobi coefficients.
\end{proof}

We make two remarks on Theorem \ref{J-lift}.

The modular forms $\nabla_3^2$ and $\nabla_2^2$ coincide with
generators of the graded rings of Siegel modular forms for
$\Gamma_0^{(2)}(2)$ and $\Gamma_0^{(2)}(3)$ (see \cite{Ib}).
The lifting construction  gives us an universal approach
to the generators. Moreover we obtain more fundamental functions
like $\nabla_2$ or $\nabla_{3/2}$ which are  roots from
generators of the corresponding graded rings.
We give the relations between dd-modular forms
and the generators proposed in the papers of Ibukiyama.
First, we have
$$
\nabla_3(Z)^2=\Lift\bigl(\eta(\tau)^2\eta(2\tau)^8\,\vartheta(\tau,z)^2\bigr)
=K(Z)\in S_6(\Gamma_0(2))
$$
where
$$
K(Z)=\frac{1}{4096}
\bigl(\theta_{0100}(Z)\theta_{0110}(Z)\theta_{1000}(Z)\theta_{1001}(Z)
\theta_{1100}(Z)\theta_{1111}(Z)\bigr)^2
$$
and
$$
\nabla_2(Z)^2=
\Lift\bigl(\eta(3\tau)^6\vartheta(\tau,z)^2\bigl)=\frac{1}{24}\,\Theta_4(Z)
\in S_4(\Gamma_0(3)).
$$
$\theta_{abcd}$ denotes  the  Siegel theta-series with characteristic $(abcd)$
of level $2$ and $\Theta_4$ is a theta-series with a spherical function.
Using the dd-function $Q_1$ we can construct the generators of the graded ring
of the modular forms with respect to the congruence subgroup
$\Gamma_2(2)$. The details will be published in a separate paper.

The second remark is related to  differential equations.
In  \cite{CYY}  it was proved that the monodromy group of
Picard--Fuchs  equations  associated with one parameter families
of Calabi--Yau threefolds is a subgroup of certain congruence subgroup
$\Gamma(d_1,d_2)$ in $\Sp_2(\ZZ)$ where $d_2$ is a divisor of $d_1$.
Therefore one can put a question on Siegel modular forms with respect to this group.
In order to construct such modular forms we can use Theorem \ref{J-lift}
because this subgroup  is the integral part of the intersection
of two modular groups considered in this theorem
$$
\Gamma(d_1,d_2)=\Gamma_{t_1}(q_1)\cap \Gamma_{t_2}(q_2)\cap \Sp_2(\ZZ)
$$
where we have the following relations for the least common multiples
$d_1=[t_1,t_2]$ and $d_2=[q_1,q_2]$.
If $t_1$ and $t_2$ are coprime we do not need to make the intersection with $\Sp_2(\ZZ)$.
In particular, $\Gamma(d_1,d_2)=\Gamma_0^{(2)}(d_2)\cap \Gamma_{d_1}$.
According to Theorem \ref{J-lift} for any
$f_i\in  J_{k_i,t_i}(\Gamma_0(q_i))$ ($i=1,\,2$)
the product ${\rm Lift}(f_1)\cdot {\rm Lift}(f_2)$ is a modular form of weight
$k_1+k_2$ with respect to $\Gamma(d_1,d_2)$.
(Under some conditions on $t_i$ and $q_i$ one can consider
Jacobi forms with some characters.)

To finish the proof of Theorem \ref{dd-class} we have to construct
a square root from the $\Gamma_0^{(2)}(4)$-modular form $F_3$
(see (\ref{F3})).
For this aim  we consider the Borcherds automorphic products.

\section{Borcherds products for $\Gamma_t(N)$}

In this section we consider Borcherds automorphic products
related to the Jacobi forms of weight $0$ with respect to the congruence subgroup
$\Gamma_0(N)$. This construction gives us the dd-modular forms of
\S 2 as automorphic products.
In particular, we construct the last dd-modular form
$\nabla_{3/2}$ with respect to $\Gamma_0^{(2)}(4)$.
In \cite{B1} the language of the orthogonal groups  and the  vector valued automorphic forms
was used.
The Jacobi forms are very useful in the framework of
Siegel modular forms because we have many methods to construct Jacobi forms
of weight $0$.
The case of the symplectic paramodular group $\Gamma_t$
was considered in \cite{GN1}--\cite{GN2}.
A similar result one can obtain for the congruence subgroups.
Some examples of Borcherds automorphic products for
$\Gamma_0^{(2)}(N)<\Sp_2(\ZZ)$  in terms of Jacobi forms
were constructed in \cite{AI}  but they could not prove that the construction works
for arbitrary Jacobi forms (see Lemma \ref{D1} below and the remark before it).
In this section we construct the automorphic products for
the subgroups $\Gamma_t(N)$ of the paramodular groups $\Gamma_t$ for any $t$ and $N$.

First we recall some well known facts about
the Hecke congruence subgroup  $\Gamma_0(N)$ (see \cite{Sh}, \cite{Mi}).
The number of non-equivalent cusps of $\Gamma_0(N)$
is equal to
$\displaystyle \sum_{\substack {e\vert N,\, e>0}} \varphi((e,{\frac{N}e}))$
where $\varphi$ is the Euler's function and $(a,b)$
is the greatest common divisor of $a$ and $b$.
We denote by ${\mathcal P}$ the set of cusps
$$
{\mathcal P}=\displaystyle \left\{{\frac{f}e},
\   e\vert N,\ e\geq 1,\   f\mod (e,{\frac{N}e}),\   (e,f)=1\right\}.
$$
To each cusp ${f}/{e}\in \mathcal P$ of $\Gamma_0(N)$,
we associate a matrix
$$
\frac{f}{e} \mapsto M_{f/e}=
\left(\begin{array}{cc} f & * \\ e & *\end{array}\right)\in \SL(2,\ZZ),
\quad M_{f/e}\langle \infty\rangle=f/e.
$$
Let $h_e={N}/(e^2,N)$ be the width of the cusp ${f}/e \in \mathcal P$.
The sum of the widths is
$N\cdot \prod_{p\,|N}(1+p^{-1})$ ($p$ is prime)
which is the index of $\Gamma_0(N)$ in $\SL(2,\ZZ)$.
We also put $N_e={\frac{N}e}$. In order to  construct the  dd-modular forms
we need two particular cases when $N=p$ or $p^2$.
\smallskip

\noindent
{\bf Example.} $\Gamma_0(p)$ and $\Gamma_0(p^2)$.

\noindent
i) If $N=p$,  $p$ prime, then there are  two cusps:
$\frac{1}{p}$ which is $\Gamma_0(p)$-equivalent to $\infty$ and $0$
of width $1$ and $p$ respectively.
\smallskip

\noindent
ii) If $N=p^2$, $p$ prime , there are  $(p+1)$ cusps:
$\frac{1}{p^2}$ which is $\Gamma_0(p^2)$-equivalent to $\infty$,
$0$ and $\left\{\frac{f}{p},\ 1\leq f\leq p-1\right\}$
of width $1$, $p^2$ and $1$ respectively.
\smallskip

As we mentioned above our datum for  the automorphic Borcherds product
for the congruence subgroup $\Gamma_t(N)<\Gamma_t$ is a nearly holomorphic Jacobi
form of weight $0$  and index $t$ with respect to $\Gamma_0(N)$ (see \S 2).
The character of this form is trivial.
In the Borcherds automorphic products  \cite{B1}
vector valued modular forms were used.
In the case of a Jacobi modular form with respect to the  congruence subgroup
$\Gamma_0(N)$ one has to use its Fourier coefficients
at all cusps of $\Gamma_0(N)$
(see \cite[Examples 2.2 and 2.3]{B1}).
In order to realize this one can  use the complete Hecke operator $T_N(m)$
for $\Gamma_0(N)$ which contains more classes than the operator $T_-^{(N)}(m)$
defined in (\ref{Tm}) if $(m,N)\ne 1$.
The operator $T_N(m)$ was introduced in \cite{He}  and it was used in  \cite{AI}.
For $m\in \NN^*$, we set
$$
M_N(m)=\left\{M=\left(\begin{array}{cc} a & b \\ cN & d\end{array}\right)
\in\mathcal M_2(\ZZ)
\,\,
|\,{\text{det}}(M)=m \right\}.
$$
Similar to (\ref{Tm})  we can consider the Hecke operator with respect
to the parabolic
subgroup of $\Gamma_t(N)$ acting on the modular forms
$\tilde \phi(Z)=\phi(\tau,z)e^{2\pi i t\omega}$. This gives us  for any
$\phi\in J_{0,t}^{nh}(\Gamma_0(N))$ the Hecke operator
\begin{equation}\label{TNm}
\hspace{-0.2cm}
\phi\big\vert_{0,t}T_N(m) (\tau,z)
=m^{-1}\hspace{-0.7cm}
\sum_{\substack
{\scriptscriptstyle
\left(\smallmatrix a&b\\c&d\endsmallmatrix\right)
\in \Gamma_0(N)\backslash M_N(m)}}\hspace{-0.4cm}
e^{-2i\pi mt\,\tfrac{cz^2}{c\tau+d}}\,\phi\,(\frac{a\tau+b}{c\tau+d}\,,
\frac{mz}{c\tau+d}).
\end{equation}
Then  $\phi\big\vert_{0,t}T_N(m)\in J_{0,mt}^{nh}(\Gamma_0(N))$.
This operator transfers the weak (holomorphic) Jacobi forms into weak
(holomorphic) Jacobi forms.

We can  write the  Fourier expansion of  $\phi\in J_{0,t}^{nh}(\Gamma_0(N))$
(see (\ref{Fexp})) at the  corresponding cusp  ${f}/{e}$ using $M_{f/e}$
$$
\bigl(\phi\big\vert_{0,t}M_{{f}/{e}}\bigr)(\tau,z)
=\sum_{ n\in \ZZ/h_e}
\sum_{l\in \ZZ}c_{{f}/{e}}(n,l)q^nr^l.
$$
We note that $c_{1/N}(n,l)$ is the Fourier coefficient of $\phi$
at infinity.
For  a weak Jacobi form we have $n\ge 0$ if $c_{{f}/{e}}(n,l)\ne 0$.
\begin {theorem}\label{product}
Let $\phi\in J_{0,t}^{nh}(\Gamma_0(N))$.
Assume that for all cusps of $\Gamma_0(N)$ we have that
$\frac{h_e}{N_e}c_{{f}/{e}}(n,l)\in \ZZ$ if  $4nmt-l^2\leq 0$.
Then the  product
$$
B_\phi(Z)=q^{A}r^{B}s^{C}\prod_{f/e \in \mathcal P}\
\prod_{\substack{n,l, m\in \ZZ \vspace {1pt} \\ (n,l,m)>0}}
\bigl(1-(q^nr^ls^{tm})^{N_e}\bigr)^{\frac{h_e}{N_e}c_{{f}/{e}}(nm,l)},
$$
where
$(n,l,m)>0$ means  that if $m>0$, then $n\in \ZZ$ and $l\in \ZZ$,
if $m=0$ and $n>0$, then  $l\in \ZZ$,
if $m=n=0$, then $l<0$, and
$$
A=\frac{1}{24}\sum_{\substack{{f}/{e}\,\in \mathcal P \vspace{0.5pt}\\l\in \ZZ }}
\hspace{-3pt}h_ec_{{f}/{e}}(0,l),\ \
B=\frac{1}{2}\sum_{\substack{{f}/{e}\,\in \mathcal P \vspace{0.5pt}\\l\in \ZZ,\, l>0 }}
\hspace{-3pt}lh_e c_{{f}/{e}}(0,l), \ \
C=\frac{1}{4}\sum_{\substack{{f}/{e}\,\in \mathcal P \vspace{0.5pt} \\l\in \ZZ  }}
\hspace{-3pt} l^2h_e c_{{f}/{e}}(0,l),
$$
defines a meromorphic modular  form of weight
$$
k=\frac{1}{2}\sum_{{f}/{e}\in \mathcal P}\frac{h_e}{N_e}\,c_{{f}/{e}}(0,0)
$$
with respect to $\Gamma_t(N)^+$ with a character
(or a multiplier system) $\chi$. In particular
$$
\frac{ B_\phi(V_t\langle Z \rangle)}{B_\phi(Z)}
=(-1)^{D_0}\quad
\text{with}\quad
D_0
=\sum_{{f}/{e}\in \mathcal P}\sum_{l\in \ZZ ,\ n<0}
\frac{h_e}{N_e}\sigma_0(-n)c_{{f}/{e}}(n,0)
$$
where
$V_t\langle Z \rangle
=V_t\langle
\left( \smallmatrix \tau&z \\z& \omega\endsmallmatrix\right)\rangle
=\left(\smallmatrix t\omega&z \\z& {\tau}/{t}\endsmallmatrix\right)$
and
$\sigma_j(m)=\sum_{d|m} d^j$.
The poles and zeros of $B_\phi$ lie on the rational
quadratic divisors defined by the Fourier coefficients
$c_{f/e}(n,l)$ with   $4nmt-l^2<0$.
In particular $B_\phi$ is holomorphic if all  such coefficients are positive.
The character $\chi$ is induces by the following relations
$$
\chi(\widetilde{M})
=\prod_{{f}/{e}\in \mathcal P}
\bigl(v_{\eta}^{(N_e)}(M)\bigr)^
{\frac{h_e}{N_e}\sum_{l\in \ZZ}c_{{f}/{e}}(0,l)}
$$
for $M\in\Gamma_0(N_e)$ where
$v_{\eta}^{(N_e)}(M)=v_{\eta}(\alpha_{e}M\alpha_{e}^{-1})$
with $\alpha_{e}=\left(\smallmatrix N_e&0\\0&1\endsmallmatrix\right)$ and
$$
\chi([\lambda, \mu;0])
=\prod_{f/e\in \mathcal P,\ l>0}
v_{H,N_e}^{\frac{h_e}{N_e}\,lc_{{f}/{e}}(0,l)}([\lambda, \mu;0])
$$
where
$v_{H,N_e}([\lambda, \mu;0])=(-1)^{\lambda+N_e\mu+N_e\lambda\mu}$
for $\lambda,\ \mu\in \ZZ$ and
for all $\kappa\in \ZZ$
$\chi([0,0;\frac{\kappa}{t}])=e^{2i\pi C{\kappa}/{t}}$.
\end{theorem}
\begin{proof}
The paramodular group $\Gamma_t$ can be realized as the stable orthogonal group
of the lattice $2U\oplus \langle -2t\rangle$ of signature $(2,3)$ where
$U=\left(\smallmatrix 0&1\\1&0\endsmallmatrix\right)$
is the hyperbolic plane (see \cite{G1}, \cite{GH2}).
Using the similar arguments we  can realize $\Gamma_t(N)$ as a subgroup
of the orthogonal group of the lattice $U\oplus U(N)\oplus \langle -2t\rangle$
where $U(N)=\left(\smallmatrix 0&N\\N&0\endsmallmatrix\right)$.
The product of the theorem is a specialization of the Borcherds automorphic product
considered in \cite[Theorem 13.3]{B1}.
It converges if $Y=\Im Z$ lies in a Weyl chamber determined by the action
of $\Gamma_0(N)$ on $Y>0$ with $\det(Y)>C$ for a sufficiently large $C$.
The product can be extended to a meromorphic function on $\HH_2$ whose poles
and zeros lie on rational quadratic divisor of $\HH_2$.
We define below the invariants (the modular  group, the weight, the character,
the first and the second Fourier-Jacobi coefficients)
of this modular form in terms of the Fourier coefficients of the lifted Jacobi form
of weight $0$ using a representation similar to \cite{GN1}--\cite{GN2}.

We have the following decomposition (see \cite{He})
\begin{equation}\label{Tmc}
\Gamma_0(N)\backslash M_N(m)=
\end{equation}
$$
\bigsqcup_{{f}/{e}\in \mathcal P}
\left\{M_{{f}/{e}}
\left(\begin{array}{cc} a & b \\ 0 &d \end{array}\right)\ |\
ad=m,\ \  ae\equiv 0 \bmod N,\ \
 b\bmod h_ed  \right\}.
$$
For $\phi\in J_{0,t}^{nh}(\Gamma_0(N))$ and
$Z=\left(\smallmatrix \tau&z\\z&\omega\endsmallmatrix\right)\in \HH_2$
we set
\begin{equation}\label{L-phi}
L_\phi(Z)=\sum_{m=1}^{\infty} \phi\vert_{0,t}T_N(m)
(\tau,z)e^{2i\pi tm \omega}.
\end{equation}
Using the decomposition of  $\Gamma_0(N)\backslash M_N(m)$
and the formula for the action of $T_N(m)$
(see the proof of Theorem \ref{J-lift})
we have (whenever the product converges):
$$
\text{Exp}(-L_\phi(Z))
=\prod_{{f}/{e}\in \mathcal P}
\prod_{\substack{m\geqslant 1 \vspace{0.5\jot} \\ n,\,l\in \ZZ}}
\bigl(1-(q^nr^ls^{tm})^{N_e}\bigr)^{\frac{h_e}{N_e}c_{{f}/{e}}(nm,l)}.
$$
This product is invariant with respect to the action of the Jacobi group.
We introduce one more factor(the ``zeroth" Hecke operator or the  Hodge
correction in the geometric terms of \cite{G4})
\begin{equation}\label{T0}
\hspace{-0.3truecm}T_\phi^{(0)}(Z)=\prod_{{f}/{e}\in \mathcal P}
\eta(N_e\tau)^{\frac{h_e}{N_e}c_{{f}/{e}}(0,0)}
\prod_{l>0}\biggl(\frac{\vartheta(N_e\tau,N_e lz)}{\eta(N_e\tau)}\,
e^{i\pi N_el^2\omega}\biggr)^{\frac{h_e}{N_e}c_{{f}/{e}}(0,l)}.
\end{equation}
So as in \cite[(2.7)]{GN2}  we obtain that
\begin{equation}\label{Bphi}
B_\phi(Z)=T_\phi^{(0)}(Z)\cdot\text{Exp}(-L_\phi(Z)).
\end{equation}
The additional term    $T_\phi^{(0)}(Z)$ is a nearly holomorphic
Jacobi form of weight $k$ indicated in the theorem
and of index $C\in \NN/2$ with respect to $\Gamma_0(N)$.
This is {\it the first Fourier-Jacobi coefficient} of the automorphic product $B_\phi$.
(It might be that this is a Jacobi form of index zero, i.e., an automorphic form in $\tau$.)
The Jacobi form is a modular form with respect to the parabolic subgroup
$\Gamma_t^{\infty}(N)$.
Like in the proof of Theorem \ref{J-lift}  we  use that
$\Gamma_t(N)^+=\langle \Gamma_t^{\infty}(N), V_t\rangle$.
We have to analyze  the behavior of   $B_\phi$ under  $V_t$-action.
Like in \cite{GN2} a straightforward calculation shows that
$$
\frac{ B_\phi(V_t\langle Z \rangle)}{B_\phi(Z)}
=(-1)^{D_0}(q^{1/t}s^{-1})^{tD_1+C-tA}
$$
where $D_0$ is given in the theorem and
$$
D_1
=\sum_{{f}/{e}\in \mathcal P}\ \ \sum_{l\in \ZZ,\  n<0}
h_e\sigma_1(-n)c_{{f}/{e}}(n,l).
$$
We note that in $\cite{AI}$ the approach of \cite{GN1}--\cite{GN2} was also  used.
But in \cite{AI} it was not proved  that $tD_1+C-tA=0$.
To show that the automorphic product is  $V_t$-invariant
we  prove  Lemma \ref{D1} (see below) similar to \cite[Lemma 2.2]{GN2}.
We note also that
the automorphic product of the theorem is  defined at the ``standard"
$0$-dimensional cusp $\infty$ of $\Gamma_t(N)$.
If $N=1$ then the $\Gamma_t$-orbit of any rational quadratic
divisor (a Humbert modular surface) has a representative containing $\infty$
(see \cite{GH2} and \cite{GN2}).
If  $N>1$ then there are more orbits. Not all of them have a non-trivial
intersection with infinity. Therefore the arguments in the construction
of some examples of the automorphic products in \cite{AI} are not complete.
One has to use \cite[Theorem 13.3]{B1} in the proof.
We add that as we mentioned in the beginning of this proof
we do not agree with \cite[page 262]{AI} that ``$\Gamma_0^{(2)}(N)$
is not an automorphism group of a lattice."

\begin{lemma}\label{D1}
For any  $\displaystyle \phi\in J_{0,t}^{nh}(\Gamma_0(N))$
we have $\displaystyle tD_1+C-tA=0$.
\end{lemma}
\begin{proof}
We give a proof based on the method of the automorphic correction proposed in \cite{G4}
 which is more simple than the proof of \cite[Lemma 2.2]{GN2}.
For any $\displaystyle \phi\in J_{k,t}^{nh}(\SL_2(\ZZ))$,
we consider the following automorphic correction of $\phi$:
$$
\Psi(\tau,z)=e^{-8\pi^2tG_2(\tau)z^2}\phi(\tau,z) \quad
\text{where }\ G_2(\tau)=-\frac{1}{24}+\sum_{n\geqslant 1}\sigma_1(n)q^n
$$
is the quasi-modular Eisenstein series of weight $2$.
The corrected form satisfies the functional equation
$$
\Psi(\frac{a\tau+b}{c\tau+d},\frac{z}{c\tau+d})=(c\tau+d)^k\Psi(\tau,z),
\qquad
\forall\ \begin{pmatrix} a & b \\ c & d\end{pmatrix}\in \SL_2(\ZZ).
$$
We consider the  Taylor expansion of $\Psi$ around $z=0$
$$
\Psi(\tau,z)=\sum_{\nu\geqslant 0}f_{\nu}(\tau)z^{\nu}.
$$
The  Taylor  coefficient $f_\nu\in M^{(mer)}_{k+\nu}(\SL_2(\ZZ))$
are modular forms with a possible  pole of finite order at the cusp.
If
$
\phi(\tau,z)=\sum_{\substack {n\in\ZZ}}\sum_{\substack {l\in \ZZ}}c(n,l)q^nr^l
$ is in $ J_{0,t}^{nh}(\SL_2(\ZZ))$,
then
$$
f_2(\tau)
=\frac{\partial^2\phi}{\partial z^2}(\tau,0)-16\pi^2tG_2(\tau)\phi(\tau,0)
\in M^{(mer)}_{2}(SL_2(\ZZ)).
$$
But the constant term of any nearly holomorphic modular form of weight two is zero
(see \cite[Lemma 9.2]{B2}). Therefore
\begin{equation}\label{=0}
t\sum_{l\in\ZZ}c(0,l)-24t\sum_{\substack {n<0 \\ l\in\ZZ}}
\sigma_1(-n)c(n,l)-6\sum_{l\in\ZZ}l^2c(0,l)=0.
\end{equation}
For a Jacobi form with respect to a congruence subgroup we use the trace operator.
Let $\phi\in J_{0,t}^{nh}(\Gamma_0(N))$. We set
$$
\psi={\rm Tr}_{\SL_2(\ZZ)}(\phi)=
\sum_{\gamma\in \Gamma_0(N)\backslash \SL_2(\ZZ)}
\phi\big\vert _{0,t}\gamma \in J_{0,t}^{nh}(\SL_2(\ZZ)).
$$
This Jacobi form has the following Fourier expansion
$$
\psi(\tau,z)=\sum_{\substack {n\in\ZZ}}\sum_{\substack {l\in \ZZ}}c(n,l)q^nr^l
\quad\text{where}\quad
c(n,l)=\sum_{{f}/{e}\in \mathcal P}h_e c_{{f}/{e}}(n,l).
$$
The last expression is obtained by noticing that
$$
\SL_2(\ZZ)
=\bigsqcup_{{f}/{e}\in \mathcal P,\  0\leqslant a\leqslant h_e-1}
\Gamma_0(N)M_{{f}/{e}}
\begin{pmatrix} 1&a\\0&1\end{pmatrix}.
$$
The claim of the lemma follows from (\ref{=0}).
\end{proof}
The formula for the character  of $\Gamma_t(N)^{+}=\langle \Gamma_t(N),V_t\rangle$
follows directly form
the calculation with $\eta$- and $\vartheta$-factors
in $T_\phi^{(0)}$ (see (\ref{T0})) which is  the first Fourier-Jacobi
coefficient of $B_\phi$. More exactly
the $\SL_2$-part of the character (or the multiplier system) of this Jacobi form
is equal to the character of the $\eta$-product
$$
\prod_{{f}/{e}\in \mathcal P}
\eta(N_e\tau)^{\frac{h_e}{N_e}\sum_{l}c_{{f}/{e}}(0,l)}
$$
which is a $\Gamma_0(N)$-modular form because $N_e$ is a divisor of $N$.
Its character is the character  $\chi(\widetilde M)$ of  the theorem.
The Heisenberg  part of the character of the Jacobi form
$\vartheta(M\tau,Mz)$ of index $M/2$ is equal to
$$
v_{H,M}([\lambda, \mu;0])=v_H([\lambda, M\mu;0])=
(-1)^{\lambda+M\mu+M\lambda\mu}.
$$
It gives us the Heisenberg part of the character.
We note that the second Fourier-Jacobi coefficient is equal to
$T_\phi^{(0)}\cdot \widetilde \phi$.
We note that if a Siegel modular form $F$ is a Borcherds automorphic product
$B_\phi$ we can find $\phi$ taking the quotient of the first two non-zero
Fourier--Jacobi coefficients of $F$.
\end{proof}

In order to obtain the Borcherds products for the dd-modular forms
we propose a method of construction of  weak Jacobi forms  of weight $0$
for $\Gamma_0(N)$ using the Jacobi theta-series with characteristics (see \cite{Mu}).
Let  $N\in \NN$ and $(a,b)\in \ZZ^2$.
We call the theta-series of level $N$ with characteristic $(a,b)$, the series
\begin{equation}\label{thetacar}
\vartheta_{a,b}^{(N)}(\tau,z)
=\sum_{n\in \ZZ}e^{i\pi(n+\frac{a}{N})^2\tau+2i\pi(n+\frac{a}{N})(z+\frac{b}{N})}.
\end{equation}
This is a holomorphic function on $\HH_1 \times\CC$.
Among these series, there is a special one for $(a,b)=(0,0)$
$$
\vartheta_{00}(\tau,z)=\sum_{n\in \ZZ}e^{i\pi n^2\tau+2i\pi nz}
=\prod_{n\geqslant1}(1-q^n)(1+q^{\frac{2n-1}{2}}r)(1+q^{\frac{2n-1}{2}}r^{-1}).
$$
All the theta-series with characteristics can be expressed by the mean of $\vartheta_{00}$
$$
\vartheta_{a,b}^{(N)}(\tau,z)
=e^{2i\pi\frac{ab}{N^2}}q^{\frac{a^2}{2N^2}}r^{\frac{a}{N}}
\vartheta_{00}(\tau,z+\frac{a}{N}\tau+\frac{b}{N}).
$$
We also have for any $(a',b')\in \ZZ^2$
$$
\vartheta_{a+a'N,b+b'N}^{(N)}(\tau,z)=e^{2i\pi \frac{ ab'}{N}} \vartheta_{a,b}^{(N)}(\tau,z).
$$
The last formula allows us to take the characteristics $(a,b)$ modulo $N$.

To construct Jacobi forms of weight $0$ we consider quotients of theta-series.
We put
$$
\xi_{a,b}^{(N)}(\tau,z)=\frac{\vartheta_{a,b}^{(N)}(\tau,z)}{\vartheta_{a,b}^{(N)}(\tau,0)}.
$$
This function is holomorphic on $\mathbb{H}_1\times\CC$
for any $(a,b)$ if  $N$  is odd.
For $N$ even, as $\vartheta_{00}(\tau,\frac{\tau}{2}+\frac{1}{2})=0$,
we cannot make the quotient $\xi_{\frac{N}{2},\frac{N}{2}}^{(N)}$.
When we write $\xi_{a,b}^{(N)}$ for even
$N$ then we assume that  $(a,b)\neq(\frac{N}{2},\frac{N}{2})$.
In fact $\xi_{a,b}^{(N)}$ is a  weak Jacobi forms of weight $0$ and index $1/2$
with respect to  the principal congruence subgroup of level $N$
(see (\ref{J-slash})):
\begin{equation}\label{modN}
\xi_{a,b}^{(N)}\vert_{0,\frac{1}{2}}M=\xi_{a,b}^{(N)},
\qquad M\in\Gamma(N).
\end{equation}
More exactly we have the following functional equations with respect
to the generators of the full Jacobi group (see \S 2)
$$
\xi_{a,b}^{(N)}\vert_{0,\frac{1}{2}}[\lambda,\mu;0]
=e^{2i\pi\frac{a}{N}\mu}e^{-2i\pi\frac{b}{N}\lambda}\xi_{a,b}^{(N)},
\qquad (\lambda,\mu)\in\ZZ^2,
$$
$$
\xi_{a,b}^{(N)}\vert_{0,\frac{1}{2}}S=\xi_{b,\overline{-a}}^{(N)}\, ,
\qquad
S=\left(\begin{smallmatrix}
0&-1\\1&\phantom{-}0\end{smallmatrix}\right)
$$
where $\overline{-a}$ is the unique representant of $-a$ modulo $N$ such that
$\overline{-a}\in\left\{0,...,N-1\right\}$,
$$
\xi_{a,b}^{(2N)}\vert_{0,\frac{1}{2}}T=\xi_{a,\overline{a+b+N}}^{(2N)}\, ,
\qquad
T=\left(\begin{smallmatrix}
1&1\\0&1\end{smallmatrix}\right)
$$
where $\overline{a+b+N}$ is the unique representant of $a+b+N$ modulo $2N$
such that
$\overline{a+b+N}\in \left\{0,...,2N-1\right\}$
and
$$
\xi_{a,b}^{(2N'+1)}\bigg\vert_{0,\frac{1}{2}}T
=\xi_{2a,\overline{2(a+b+N')+1}}^{(4N'+2)}
$$
where $\overline{2(a+b+N')+1}$ is the unique representant of $2(a+b+N')+1$
modulo $4N'+2$
such that $\overline{2(a+b+N')+1}\in \left\{0,...,4N'+1\right\}$.
These formulae lead us to construct weak Jacobi forms for $\Gamma_0(N)$
in the following way:

\noindent
(i) we consider the quotient group $G=\Gamma(N)\backslash\Gamma_0(N)$
if $N$ is even
or $G=\Gamma(2N)\backslash\Gamma_0(N)$ if $N$ is odd since according
to the $T$-transformation formula we have to double the level;

\noindent
(ii) we compute the orbits of $\xi_{a,b}^{(N)}$ under $G$;

\noindent
(iii) in a fixed  orbit of $\xi_{a,b}^{(N)}$, we take some powers of elements
or products of them
in order to obtain the trivial character of the Jacobi group.
\smallskip

\noindent
In this paper we only construct the Jacobi forms of weight $0$
which generate dd-modular forms.
We are planing  to obtain results similar to \cite{G4} about the structure of the graded rings
of weak Jacobi forms with respect to $\Gamma_0(N)$ for small $N$
in a separate publication.
\smallskip

\noindent
{\bf Examples 3.3} {$\mathbf N=\mathbf 2$}.
We have $G=\left\{I_2, T\right\}$ (the group of order two)
and the orbit $O_G(\xi_{1,0}^{(2)})$ contains the only element
$\xi_{1,0}^{(2)}$.
The formula for the $[\mu,\nu;0]$-action  implies that
$ \xi_{1,0}^{(2)}$ has a character of order two. Therefore
$$
(\xi_{1,0}^{(2)})^2\in J_{0,1}^{w}(\Gamma_0(2)).
$$
{$\mathbf N=\mathbf 3$}.
In this case $G$ is non abelian group of order 36.
It contains the set
$\Sigma=\{\pm T^k, \pm ST^3ST^k,0\leqslant k \leqslant 5\}$.
Therefore
$O_G(\xi_{3,1}^{(6)})\supseteq O_\Sigma(\xi_{3,1}^{(6)})
=\{\xi_{3,1}^{(6)},\xi_{3,5}^{(6)}\}$
and using the standard generators of $\Gamma_0(3)$ we have equality.
Therefore
$$
\xi_{3,1}^{(6)}\cdot \xi_{3,5}^{(6)}\in J_{0,1}^{w}(\Gamma_0(3)).
$$
{$\mathbf N=\mathbf 4$}.
We have that $G=\{I_2,T,T^2,T^3,-I_2,TST^4S,TST^4ST,TST^4ST^2\}$
is the group of order 8 isomorphic to $\ZZ/2\ZZ\times \ZZ/4\ZZ$.
We see that $O_G(\xi_{0,1}^{(4)})=\{\xi_{0,1}^{(4)},\xi_{0,3}^{(4)}\}$
and $O_G(\xi_{2,1}^{(4)})=\{\xi_{2,1}^{(4)},\xi_{2,3}^{(4)}\}$.
Therefore
$$
\xi_{0,1}^{(4)}\cdot \xi_{0,3}^{(4)},\
\text{ and }\ \xi_{2,1}^{(4)}\cdot \xi_{2,3}^{(4)}
\in J_{0,1}^{w}(\Gamma_0(4)).
$$

{\bf The dd-modular forms as Borcherds products.}
Now we can finish the proof of Theorem \ref{dd-class} and to construct
the last dd-modular form $\nabla_{3/2}$ of weight $3/2$ with respect to $\Gamma_0^{(2)}(4)$.
Moreover we give the Borcherds automorphic product for all new  Siegel
dd-modular forms with respect to the congruence subgroups.
(The Borcherds products of the  dd-modular forms for the full paramodular group
were found in  \cite{GN1}--\cite{GN2}.)

We start with $N=2$.
Let
$$
\phi_2(\tau,z)=4(\xi_{1,0}^{(2)})^2(\tau,z)\in J_{0,1}^{w}(\Gamma_0(2)).
$$
There are two cusps and $\phi_2$ has the two Fourier expansions with
integral Fourier coefficients
$$
\phi_2(\tau,z)=(r^{-1}+2+r)+2(r^{-2}-2+r^2)q+\dots
=\sum_{n\in \NN,\,l\in \ZZ}c(n,l)q^nr^l,
$$
$$
(\phi_2\vert_{0,1}S)(\tau,z)=4-8(r^{-1}-2+r)q^{\frac{1}{2}}+\dots
=\sum_{n\in \frac{1}{2}\NN,\,l\in \ZZ}c_S(n,l)q^nr^l .
$$
The only orbit of the Fourier coefficients with negative hyperbolic norm $4nt-l^2$
of its index is $c(0,1)=1$.
Then applying Theorem \ref{product} to $\phi_2$, we obtain
$$
B_{\phi_2}(Z)=q^{\frac{1}{2}}r^{\frac{1}{2}}s^{\frac{1}{2}}
\prod_{(n,l,m)>0}(1-q^nr^ls^m)^{c(nm,l)}(1-q^{2n}r^{2l}s^{2m})^{c_S(nm,l)}
$$
$$
=\eta(\tau)\eta(2\tau)^4\vartheta(\tau,z)e^{i\pi\omega}\cdot{\text{Exp}}(-L_{\phi_2})(Z).
$$
This is a holomorphic Siegel modular form of weight $3$ with respect to $\Gamma_0(2)$.
According to the Koecher principle a Siegel dd-modular form is defined up to a constant.
Comparing the first Fourier coefficients we obtain
$$
\nabla_3(Z)=\Lift(\eta(\tau)\eta(2\tau)^4\vartheta(\tau,z))=B_{\phi_2}(Z).
$$

For $N=3$ we take
$$
\phi_3(\tau,z)=3(\xi_{3,1}^{(6)}\xi_{3,5}^{(6)})(\tau,z)\in J_{0,1}^{w}(\Gamma_0(3)).
$$
We again  have two Fourier expansions containing only integral Fourier coefficients
$$
\phi_3(\tau,z)=(r^{-1}+1+r)+(r^{-2}-r^{-1}-r+r^2)q+\dots
=\sum_{n\in \NN,\,l\in \ZZ}c(n,l)q^nr^l,
$$
$$
(\phi_3\vert_{0,1}S)(\tau,z)=3-3(r^{-1}-2+r)q^{\frac{1}{3}}+\dots
=\sum_{n\in \frac{1}{3}\NN,\,l\in \ZZ}c_S(n,l)q^nr^l.
$$
The both Fourier expansions contain only one type of coefficients with negative norm
of its index.
This is $c(0,1)=1$.
According to Theorem \ref{product}  we obtain
$$
B_{\phi_3}(Z)=q^{\frac{1}{2}}r^{\frac{1}{2}}s^{\frac{1}{2}}
\prod_{(n,l,m)>0}(1-q^nr^ls^m)^{c(nm,l)}(1-q^{3n}r^{3l}s^{3m})^{c_S(nm,l)}
$$
and
$$
\nabla_2(Z)=\Lift(\eta(3\tau)^3\vartheta(\tau,z))=B_{\phi_3}(Z).
$$

{\bf The dd-modular form $\nabla_{3/2}$.}
The  case of $N=4$ is a little bit more difficult because there are three  different cusps.
Let
$$
\phi_4(\tau,z)=2(\xi_{2,1}^{(4)}\xi_{2,3}^{(4)})(\tau,z)\in J_{0,1}^{w}(\Gamma_0(4)).
$$
We have the following Fourier expansions
$$
\phi_4(\tau,z)=(r^{-1}+r)+(r^{-3}-r^{-1}-r^1+r^3)q^2+\dots=
\sum_{n\in \NN,\,l\in \ZZ}c(n,l)q^nr^l,
$$
$$
(\phi_4\vert_{0,1}S)(\tau,z)=2-2(r^{-1}-2+r)q^{\frac{1}{4}}+\dots
=\sum_{n\in \frac{1}{4}\NN,\,l\in \ZZ}c_S(n,l)q^nr^l,
$$
$$
(\phi_4\vert_{0,1}M)(\tau,z)=2+2(r^{-2}-2+r^2)q+\dots=
\hspace{-0.1cm}\sum_{n\in \NN,\,l\in \ZZ}c_M(n,l)q^nr^l
$$
where
$M=\left(\begin{smallmatrix} 1 & -1 \\ 2 & -1\end{smallmatrix}\right)$.
All Fourier coefficients are integral and  there exists the only type
of coefficients with negative index  norm $c(0,1)=1$.
We obtain the  Siegel modular form
$\nabla_{3/2}=B_{\phi_4}$ of weight $3/2$ for $\Gamma_0^{(2)}(4)$  given by
$$
\nabla_{3/2}(Z)=B_{\phi_4}(Z)
=\frac{\eta(2\tau)\eta(4\tau)^2}{\eta(\tau)}\,\vartheta(\tau,z)e^{i\pi\omega}
\,{\text{Exp}}(-L_{\phi_4})(Z)=
$$
$$
q^{\frac{1}{2}}r^{\frac{1}{2}} s^{\frac{1}{2}}\hspace{-0.3cm}
\prod_{(n,l,m)>0}(1-q^nr^ls^m)^{c(nm,l)}(1-q^{2n}r^{2l}s^{2m})^{\frac{1}{2}c_M(nm,l)}
(1-q^{4n}r^{4l}s^{4m})^{c_S(nm,l)}.
$$
The modular form $\nabla_{3/2}$ is the last Siegel dd-modular form which
we need in order {\it to finish the proof of Theorem \ref{dd-class}}.
Using the Koecher principle we obtain that
$$
\nabla_{3/2}(Z)^2=B_{\phi_4}^2(Z)=
\Lift\bigl(\frac{\eta(2\tau)^2\eta(4\tau)^4}{\eta(\tau)^2}\,{\vartheta(\tau,z)^2}\bigr).
$$

The  last example is the automorphic product of  the dd-modular form
$Q_1$  with  $N=2$, $t=2$ and  $k=1$.
Let
$$
\psi(\tau,z)=2\xi_{1,0}^{(2)}(\tau,2z)\in J_{0,2}^{w}(\Gamma_0(2)).
$$
We have the following Fourier expansions
$$
\psi(\tau,z)=(r^{-1}+r)+(r^{-3}-r^{-1}-r+r^3)q+
\dots =\sum_{n\in \NN,\,l\in \ZZ}c(n,l)q^nr^l,
$$
$$
(\psi\big\vert_{0,2}S)(\tau,z)=2-2(r^{-2}-2+r^2)
q^{\frac{1}{2}}
-4(r^{-2}-2+r^2)q^1-8(r^{-2}-2+r^2)q^{\frac{3}{2}}+\dots
$$
$$
=\sum_{n\in \frac{1}{2}\NN,\,l\in \ZZ}c_S(n,l)q^nr^l .
$$
Then applying Theorem \ref{product}, we obtain
$$
B_{\psi}(Z)=q^{\frac{1}{4}}r^{\frac{1}{2}} s^{\frac{1}{2}}
\prod_{(n,l,m)>0}(1-q^nr^ls^{2m})^{c(nm,l)}
(1-q^{2n}r^{2l}s^{4m})^{c_{S}(nm,l)}
$$
$$
=Q_1(Z)=\Lift\bigl(\frac{\eta(2\tau)^2}{\eta(\tau)}\vartheta(\tau,z)\bigr).
$$
\smallskip

{\bf A traced form of Borcherds product and reflective modular forms.}
For each dd-modular form  we have the identity between the known
(due to the Jacobi lifting) Fourier expansion and the Borcherds products.
We note that such examples are rather rare.
Below we give more examples of this type analyzing new
reflective modular forms, i.e., the modular forms with  divisor determined by
some reflections in the corresponding modular group (see \cite{GN2}--\cite{GN3}).
Every   dd-modular form is reflective.
We construct new examples  as the  quotient of dd-modular forms.
To represent the quotient of two dd-modular functions  in a better form we give
a new representation for the automorphic product in Theorem \ref{product}.
For that we rewrite the full Hecke operator $T_N(m)$ using
the summation with respect to the classes  from the same
subgroup $\Gamma_0(N_a)$ where $N_a=N/(a,N)$:
$$
T_N(m)=
\sum_{a|m} \ \ \sum_{M\in \Gamma_0(N)\setminus \Gamma_0(N_a)}\ \
\sum_{b \text{ mod } m/a}
\Gamma_0(N) M \begin{pmatrix}
a&b\\ 0& m/a \end{pmatrix}.
$$
Let us reorganize the formal Hecke sum
$
L_T=\sum_{m=1}^{\infty} m^{-1}T_N(m)
$
using the last representation. Formally we have
$$
L_T=
\sum_{e|N} \hspace{-0.1truecm}
\sum_{\substack{ a'\ge 1\vspace{0.5\jot} \\
(a', N_e)=1\vspace{0.5\jot}\\(a=ea')}}
\hspace{-0.1truecm}
\sum_{M\in \Gamma_0(N)\setminus \Gamma_0({N}_e)}
\sum_{\substack{ n\ge 1\vspace{0.5\jot}\\ (m=an)}}
(an)^{-1}\hspace{-0.2truecm}\sum_{b\text{ mod }n}
\Gamma_0(N)M
\begin{pmatrix}
a&0&b&0\\
0&an&0&0\\
0&0&n&0\\
0&0&0&1
\end{pmatrix}.
$$
We can rewrite the last class as
$$
\Gamma_0(N)M
\begin{pmatrix}
\frac{a}e&0&b&0\\
0&\frac{an}e&0&0\\
0&0&n&0\\
0&0&0&1
\end{pmatrix} \cdot
\begin{pmatrix}
e&0&0&0\\
0&e&0&0\\
0&0&1&0\\
0&0&0&1
\end{pmatrix}.
$$
Therefore we have a new representation for  (\ref{L-phi})
\begin{equation}\label{L-phi2}
L_\phi(Z)=
\sum_{e|N}
\sum_{\substack{m\ge 1\vspace{0.5\jot}}}
e^{-1}
\bigl(\widetilde\psi_{N_e}\vert_{\,0}\,T_-^{(N_e)}(m)\bigl)
(eZ)=
\sum_{e|N}e^{-1}L_{\psi_{N_e}}(eZ)
\end{equation}
where
$$
\widetilde\psi_{N_e}(Z)=\psi_{N_e}(\tau,z)e^{2\pi i t\omega}=
\text{Tr}_{\Gamma_0({N}_e)} \widetilde \phi(Z)=\hspace{-0.3truecm}
\sum_{M\in \Gamma_0(N)\setminus \Gamma_0({N}_e)}
\bigl(\widetilde \phi\,|_0 \,\widetilde M\bigr)(Z)
$$
is a Jacobi form of weight $0$ and index $t$ with respect to $\Gamma_0(N_e)$
and $T_-^{(N_e)}(m)$ is  the Hecke operator which we used
in the additive lifting in \S 2.
$$
T^{(N_e)}(m)=\sum_
{\substack{ad=m,\ (a, N_e)=1
\vspace{0.5\jot} \\b\,{\rm mod} \,d}}
\Gamma_0(N_e)\,
\begin{pmatrix} a&b\\0&d\end{pmatrix}.
$$
We consider  the Fourier expansion of the traced Jacobi form
$\text{Tr}_{\Gamma_0({N}_e)} \widetilde \phi$
at infinity
$$
\psi_{N_e}(\tau,z)
=\sum_{\substack{n,\,l\in \ZZ}}f_{N_e}(n,l)q^nr^l.
$$
We note that for $e=N$ we have $\psi_1=\phi$ and
$f_1(n,l)$ is the Fourier coefficient  of $\phi$ at infinity
denoted by $c_{1/N}(n,l)$ in Theorem \ref{product}.
As in the proof of
Theorem \ref{product} we have
$$
e^{-1}L_{\psi_{N_e}}(eZ)=
\sum_{\substack{m\ge 1 \vspace{0.5\jot}\\ n,l \in \ZZ}}\ \
\sum_{\substack{a\ge 1 \vspace{0.5\jot} \\ (a,N_e)=1}}
\frac{1}{ae}f_{{N}_e}(mn,l)\left(q^nr^ls^{tm}\right)^{ea}
$$
and
$$
\sum_{\substack{m\ge 1 \\ (m,N)=1}}\frac{x^m}{m}=
-\sum_{b|N}\frac{\mu(b)}{b}{\text{Log}}(1-x^b)
$$
where $\mu$ stands for the Moebius function.
Therefore
$$
L_\phi(Z)=-\sum_{e|N}\ \sum_{b|N_e}\
\sum_{\substack{m\ge 1\vspace{0.5\jot} \\ n,\,l \in \ZZ}}
{\text{Log}}\left(1-\bigl(q^nr^ls^{tm}\bigr)^{be}\right)^
{\mu(b)\frac{f_{{N}_e}(mn,l)}{be}}.
$$
The advantage of this new representation of  the Borcherds product
is evident. We use in it  only the Fourier expansion of the traced Jacobi forms
$\phi_{N_e}$ at infinity.
For the group $\Gamma_0^{(2)}(p)$ this expression contains only two functions
and one of them is  well known.
$$
\text{Tr}_{\SL_2(\ZZ)}:  J^w_{0,1}(\Gamma_0(p)) \to
J^w_{0,1}(\SL_2(\ZZ))=\CC\phi_{0,1}
$$
where
$$
\phi_{0,1}(\tau,z)=-\frac{3}{\pi^2}
\frac{\wp(\tau,z)\vartheta(\tau,z)^2}{\eta(\tau)^6}=
\sum_{n\ge 0,\,l\in \ZZ} a(n,l)q^nr^l=
(r+10+r^{-1})+\dots
$$
is one of the main  generators of the graded ring of weak Jacobi forms
(see \cite{EZ}, \cite{G4}). We note that $\phi_{0,1}$ is the elliptic genus
of Enriques surfaces and $2\phi_{0,1}$ is the elliptic genus of $\Kthree$ surfaces.
For any  $\phi_p\in  J^{nh}_{0,t}(\Gamma_0(p))$
we have
$$
\text{Exp}(-L_{\phi_p}(Z))=
\prod_{\substack{m\ge 1 \vspace{0.5\jot} \\ n,\,l\in \ZZ }}
\bigl(1-q^nr^ls^{tm}\bigr)^{c_{\phi_p}(nm,l)}
\bigl(1-q^{pn}r^{pl}s^{pmt}\bigr)^{\frac{1}p (f(nm,l)-c_{\phi_p}(nm,l))}
$$
where $c_{\phi_p}(n,l)$ and $f(n,l)$ are the Fourier coefficients
of $\phi_p$ and $\text{Tr}_{\SL_2(\ZZ)}(\phi_p)$ at infinity.

Let us consider  $\nabla_3$ ($N=2$) and $\nabla_2$ ($N=3$).
By comparing the Fourier expansions
we conclude that
$$
\phi_{0,1}=\text{Tr}_{\SL_2(\ZZ)}\phi_2
=4(\xi_{1,0}^{(2)})^2+4(\xi_{0,1}^{(2)})^2+4(\xi_{0,0}^{(2)})^2,
$$
$$
\phi_{0,1}=\text{Tr}_{\SL_2(\ZZ)}\phi_3
=3(\xi_{3,1}^{(6)}\xi_{3,5}^{(6)})+3(\xi_{1,3}^{(6)}\xi_{5,3}^{(6)})
+3(\xi_{1,1}^{(6)}\xi_{5,5}^{(6)})+3(\xi_{1,5}^{(6)}\xi_{5,1}^{(6)}).
$$
Moreover
$$
J_{0,1}^{w}(\Gamma_0(2))=\langle\phi_{0,1},\,\phi_2\rangle_{\CC}
\quad\text{and}\quad
J_{0,1}^{w}(\Gamma_0(3))=\langle\phi_{0,1},\, \phi_3\rangle_{\CC}.
$$
Therefore for $p=2$ or $3$
$$
\text{Exp}(-L_{\phi_p}(Z))=\vspace{-0.2truecm}
\prod_{\substack{m\ge 1 \vspace{0.5\jot} \\ n,\,l\in \ZZ }}
\bigl(1-q^nr^ls^{pm}\bigr)^{c_{\phi_p}(nm,l)}
\bigl(1-q^{pn}r^{pl}s^{pm}\bigr)^{\frac{1}p (a(nm,l)-c_{\phi_p}(nm,l))}
$$
where $a(n,l)$ is the Fourier coefficient of $\phi_{0,1}$.

Using this approach we can easy  calculate the product formulae
for new reflective  modular forms
of weight $2$ for $\Gamma_0(2)$, weight $3$ for $\Gamma_0(3)$,
weight $3/2$ and $7/2$ for $\Gamma_0(4)$ and  weight $1$ for $\Gamma_2(2)$:
$$
\frac{\Delta_5(2Z)}{\nabla_3(Z)},
\quad  \frac{\Delta_2(2Z)}{Q_1(Z)},
\quad \frac{\nabla_3(2Z)}{\nabla_{3/2}(Z)},
\quad \frac{\Delta_5(2Z)}{\nabla_{3/2}(Z)}
$$
and
$$
\frac{\Delta_5(Z)}{\nabla_3(Z)}, \quad \frac{\Delta_5(Z)}{\nabla_2(Z)},\quad
\frac{\Delta_2(Z)}{Q_1(Z)},
\quad \frac{\nabla_3(Z)}{\nabla_{3/2}(Z)},
\quad \frac{\Delta_5(Z)}{\nabla_{3/2}(Z)}.
$$
The dd-modular forms and all these reflective modular forms are related
to Lorentzian Kac--Moody super Lie algebras of Borcherds type.
This object  will be similar to the algebras constructed in \cite{GN1}--\cite{GN4}.
We are planning to consider  them in a  separate publication.

Using the formula $\Delta_5=B_{\phi_{0,1}}$ (see \cite[(2.16)]{GN2})
and the trace formula for $\phi_2$,
we deduce an infinite product expansion
$$
\frac{\Delta_5(Z)}{\nabla_3(Z)}=
\frac{\eta(\tau)^8}{\eta(2\tau)^4}\Prod_{\substack{m\ge 1 \vspace{0.5\jot}\\ n,\,l\in \ZZ}}
\left(\frac{1-q^nr^ls^m}{1+q^nr^ls^m}\right)^{\frac12((a(nm,l)-c_{\phi_2}(nm,l))}
$$
where $a(n,l)$ and $c_{\phi_2}(n,l)$ are respectively the Fourier coefficients
of $\varphi_{0,1}$ and $\phi_2$ at $\infty$. For $N=3$ we obtain
$$
\frac{\Delta_5(Z)}{\nabla_2(Z)}=\frac{\eta(\tau)^9}{\eta(3\tau)^3}
\Prod_{\substack{m\ge 1 \vspace{0.5\jot}\\ n,\,l\in \ZZ}}(1-q^nr^ls^m)^{b(nm,l)}
(1-q^{3n}r^{3l}s^{3m})^{-\frac13b(nm,l)}
$$
where $b(n,l)=a(nm,l)-c_{\phi_3}(nm,l)$.
The both modular forms are holomorphic because the divisor
of $\Delta_5(Z)$ is larger than the  divisor of $\nabla_3(Z)$ or $\nabla_2(Z)$.
They are non-cusp forms because the zeroth Fourier-Jacobi coefficient is non zero.

Analyzing the  examples of the reflective modular forms constructed above and
in \cite{GN2}--\cite{GN3}  we see that the first non-zero coefficient
of the Taylor expansion of a reflective form $F$ at $z=0$ is
an $\eta$-product or an $\eta$-quotient of the type considered by J. McKay
and Y. Martin (see \cite{Ma}).
We can assume that {\it every $\eta$-quotients of this type
is the  first  coefficient
of a Taylor expansion of some power of a  reflective modular form}.

The reflective modular forms in the first line above are more regular. Then we have
$$
\frac{\Delta_5(2Z)}{\nabla_3(Z)}=\widetilde\phi_{2,\frac{1}2}(Z)
\Prod_{\substack{m\ge 1 \vspace{0.5\jot}\\ n,\,l\in \ZZ}}
(1-q^{2n}r^{2l}s^{2m})^{\frac12(a(nm,l)+c_{\phi_2}(nm,l))}
(1-q^nr^ls^m)^{-c_{\phi_2}(nm,l)}
$$
where
$$
\phi_{2,\frac{1}2}(\tau,z)=
\frac{\eta(2\tau)^5}{\eta(\tau)}\frac{\vartheta(2\tau,2z)}{\vartheta(\tau,z)}
\in J_{2,\frac{1}2}(\Gamma_0(2),\chi_2)
$$
is a Jacobi cusp form of weight $2$ with a character of order $2$.
More exactly,
$\chi_2(\left(\smallmatrix a&b\\2c&d\endsmallmatrix\right))=(-1)^b$
and $\chi_2([\lambda,\mu;0])=(-1)^\lambda$.
This reflective form and its square  are the lifting
of the first Fourier--Jacobi coefficient
$$
\frac{\Delta_5(2Z)}{\nabla_3(Z)}=\text{Lift}(\phi_{2,\frac{1}2})
\in M_2(\Gamma_0^{(2)}(2),\chi_2),
$$
\begin{equation}\label{5-3}
\frac{\Delta_5(2Z)^2}{\nabla_3(Z)^2}=\text{Lift}(\phi_{2,\frac{1}2}^2)
\in M_4(\Gamma_0^{(2)}(2)).
\end{equation}
We have a  similar formula for $N=3$
\begin{equation}\label{5-2}
\frac{\Delta_5(3Z)}{\nabla_2(Z)}=\text{Lift}(\phi_{3,1})
\in M_3(\Gamma_0^{(2)}(3),\left(\frac{\det D}{3}\right))
\end{equation}
where
$$
\phi_{3,1}(\tau,z)=\eta(3\tau)^6\,\frac{\vartheta(3\tau,3z)}{\vartheta(\tau,z)}
\in J_{3,1}(\Gamma_0(3),\left(\frac{d}{3}\right)),
$$
$$
\frac{\Delta_5(3Z)}{\nabla_2(Z)}=\widetilde\phi_{3,1}(Z)
\Prod_{\substack{m\ge 1 \vspace{0.5\jot}\\ n,\,l\in \ZZ}}
(1-q^{3n}r^{3l}s^{3m})^{\frac13(2a(nm,l)+c_{\phi_3}(nm,l))}
(1-q^nr^ls^m)^{-c_{\phi_3}(nm,l)}.
$$
For $N=4$ we get two new  traced functions defined by
$\phi_4=2\xi_{2,1}^{(4)}\xi_{2,3}^{(4)}$.
They are
$$
\phi_{0,1}=\text{Tr}_{\SL_2(\ZZ)}\phi_4\quad
\text{and}\quad
\phi_2=\text{Tr}_{\Gamma_0(2)}\phi_4=\sum_{M\in \Gamma_0(4)
\setminus \Gamma_0(2)}\phi_4|_0 M.
$$
To get the second identity we take into account that
$\dim J^w_{0,1}(\Gamma_0(2))=2$.
So we are able to write the infinite product expansions for the four
reflective modular forms of type
${\Delta_5}/{\nabla_{3/2}}$  and  ${\nabla_3}/{\nabla_{3/2}}$
from our list
using the Fourier coefficients of $\phi_{0,1}$, $\phi_2$ and $\phi_4$ at infinity.

We finish with  the case $N=2$ and  $t=2$.
In order to construct $Q_1$ we used $\psi=2\xi_{1,0}^{(2)}(\tau,2z)$.
As before, we get only one new traced  function
$$
\phi_{0,2}=\text{Tr}_{\SL_2(\ZZ)}\psi=
\psi+2\xi_{0,1}^{(2)}(\tau,2z)+2\xi_{0,0}^{(2)}(\tau,2z)
$$
where $\phi_{0,2}\in J^w_{0,2}(\SL_2(\ZZ))$ is the second
generator of the graded ring of the weak  Jacobi forms of weight $0$
with integral Fourier coefficients (see \cite[(2.18)]{GN2} and \cite{G4}).
Then we get a reflective holomorphic modular form of weight $1$ with respect
to $\Gamma_2(2)<\Gamma_2$
$$
\frac{\Delta_2(Z)}{Q_1(Z)}=
\frac{\eta(\tau)^4}{\eta(2\tau)^2}
\Prod_{\substack{m\ge 1\vspace{0.5\jot} \\ n,\,l\in \ZZ}}
\left(\frac{1-q^nr^ls^{2m}}{1+q^nr^ls^{2m}}\right)^{\frac12(a_2(nm,l)-c_{\psi}(nm,l))}
$$
where $a_2(n,l)$ is  the Fourier coefficient of $\phi_{0,2}$.
In the same way we obtain that
$$
\frac{\Delta_2(2Z)}{Q_1(Z)}=\widetilde\phi_{1,\frac12}(Z)
\Prod_{\substack{m\ge 1\vspace{0.5\jot} \\ n,\,l\in \ZZ}}
(1-q^{2n}r^{2l}s^{4m})^{\frac12(a_2(nm,l)+c_{\psi}(nm,l))}
(1-q^nr^ls^{2m})^{-c_{\psi}(nm,l)}
$$
where
$$
\phi_{1,\frac12}(\tau,z)=
\eta(2\tau)\eta(\tau)\frac{\vartheta(2\tau,2z)}{\vartheta(\tau,z)}
\in J_{1,\frac{1}2}(\Gamma_0(2),\chi_4).
$$
This reflective modular form of weight one
has elementary Fourier coefficients like $Q_1$.
The character of $\phi_{1,\frac12}$ is given by the following formula
$$
\chi_4(M)=e^{\frac{2i\pi}{4}(bd+d-1)}
$$
for $M=\left(\smallmatrix a&b\\2c&d\endsmallmatrix\right)\in \Gamma_0(2)$.
Then we have
$\Gamma_1(8,4)\subset{\text{Ker}}(\chi_4)$ so $q=4$.
We also have
$$
\left(\phi_{1,\frac12}|_{\frac12}[\lambda,\mu;0]\right)(\tau,z)
=(-1)^{\mu}\phi_{1,\frac12}(\tau,z).
$$
Then we obtain that
$\frac{\Delta_2(2Z)}{Q_1(Z)}={\rm Lift} (\phi_{1,\frac12})$.
This is not a cusp form because
$
\phi_{1,\frac12}(\tau,z)
=\frac{1}{2}\vartheta_{1,0}^{(2)}(\tau,z)\vartheta_{1,0}^{(2)}(\tau,0)
$.
For $(a,8)=1$, we have
$\chi_4(\sigma_a)=\left(\frac{-4}{a}\right)$
then we deduce as for $Q_1$ that
$$
\frac{\Delta_2(2Z)}{Q_1(Z)}=\text{Lift}(\phi_{1,\frac12})(Z)
=\frac{1}{2}\sum_{N \ge 1}\sum_{\substack{n,m \in 4\NN+1 \\ l\in 2\ZZ+1 \\2nm-l^2=N^2}}
\sum_{\substack{a|(n,l,m) \\ a>0}}\left(\frac{-4}{a}\right)
q^{\frac{n}{4}}r^{\frac{l}{2}}s^{\frac{m}{2}}.
$$

\bibliographystyle{alpha}

\bigskip
\noindent
 V.~Gritsenko and F.~Cl\'ery\\
University Lille 1\\
Laboratoire Paul Painlev\'e\\
F-59655 Villeneuve d'Ascq, Cedex\\
France\\
{\tt Valery.Gritsenko@math.univ-lille1.fr}\\
{\tt Fabien.Clery@math.univ-lille1.fr}
\end{document}